\documentclass{amsart}
\usepackage{amsmath,amsthm,amssymb,amscd}

\def\({\left(}
\def\){\right)}

\def\lb\{{\left\{}
\def\rb\}{\right\}}

\newtheorem{thm}{Theorem}[section]
\newtheorem{lem}[thm]{Lemma}
\newtheorem{cor}[thm]{Corollary}

\begin{document}

\title [Enumerating orbits] {Enumerating orbits of finite classical groups on their Lie algebras}

\author{Jason Fulman}
\address{Department of Mathematics, University of Southern California, Los Angeles, CA 90089-2532, USA}
\email{fulman@usc.edu}

\author{Robert M. Guralnick}
\address{Department of Mathematics, University of Southern California, Los Angeles, CA 90089-2532, USA}
\email{guralnic@usc.edu}

\author{Richard Stong}
\address{CCR La Jolla, 4320 Westerra Ct, San Diego, CA 92121, USA}
\email{stong@ccr-lajolla.org}

\date{November 10, 2026}

\thanks{Fulman was partially supported by Simons Grant 917224. Guralnick was partially supported by a Simons Foundation Fellowship 00019819.  We thank 
Eammon O'Brien for very helpful computations and Skip Garibaldi for suggesting we point out the connection with the commuting variety. We thank ChatGPT for proofreading this paper.}

\begin{abstract} We enumerate the orbits of a finite classical group in the adjoint action on its Lie algebra, giving exact generating functions and  very good upper bounds. We treat the general linear, unitary, symplectic and orthogonal groups, in both odd and even characteristic.   In even characteristic, we get very nice formulas for the number of orbits on nilpotent classes both in the finite and algebraically closed case. We discuss the connection with commuting varieties.
\end{abstract}

\maketitle

\tableofcontents

\section{Introduction}

Let $G$ be a finite classical group.  There has been quite a lot of literature about counting the number of orbits of $G$ acting on itself by conjugation -- i.e. counting
the number of conjugacy classes of $G$.  Wall \cite{Wall} had very strong results about generating functions for this.  This was improved for the orthogonal groups
in characteristic $2$ in \cite{FG}.  It is not always easy to get bounds from these generating functions but quite good bounds were obtained in \cite{FG} and 
these bounds were used to get results about derangements and have been used extensively since, for example to study random walks on Chevalley groups and in the final solution to the Ore conjecture about commutators in finite simple groups.

In this paper, we consider a similar problem for $G$ acting on its Lie algebra.  We write down the analogous generating functions and give bounds for the number 
of orbits.  In characteristic not $2$, the answers are of the same order of magnitude as for the number of conjugacy classes (but not exactly the same).   However, in characteristic $2$ for the symplectic and  orthogonal cases, 
the answer is quite  different.   In good characteristic, using the Springer correspondence, the answers are the same for unipotent classes and nilpotent classes.  This is far from true in characteristic $2$. 

In all the cases we consider in this paper we prove  that aside from the cases of orthogonal and symplectic groups in characteristic $2$, 
the number of orbits  of the finite classical group on its Lie algebra is at least $q^n$ and at most $q^n + Cq^{n-1}$.  Here $n$ is the rank of the group
and $C$ is an explicit small constant   not depending on $n$.  

Let us describe the results in characteristic $2$. 
We first determine generating functions for the number of orbits  of nilpotent elements in the Lie algebras of the orthogonal and symplectic Lie algebras in characteristic $2$.  The answers are surprisingly nice especially given the complicated description given in \cite{LS} (see also \cite{He}).      This gives  some nice  results about the number of nilpotent orbits in characteristic $2$ for these cases both in the algebraically closed case and over finite fields (see also \cite{Xu}).  

We then deduce generating functions for the number of orbits.    This gives rise to the result that the total number of orbits is on the order of $(n+1)q^n$ for the Lie algebra
$\mathfrak{sp}(2n,q)$ and $(n+1)q^n/2$ for $\mathfrak{so}^{\epsilon}$  over a  field of size $q$ a power of $2$ (under $O$ in the last case).  There is also an asymptotic result
for $q$ fixed and the rank increasing.  

In the remaining cases, we shall see the total number of classes is close to $q^n$ (as is the case for groups ---  a result of Fulman and Guralnick \cite{FG}).  

Typically the semisimple classes are easy to describe and so the main difficulty is counting the nilpotent classes in the centralizer of a semisimple element.  The case
of semisimple classes in $\mathfrak{so}(2n,2^a)$ is surprising.  The description of these classes is fairly straightforward but  there are considerably more classes than
in odd characteristic. 
 
 We state some of our results here in characteristic $2$.  See the various sections for the results in general. The key tool in characteristic $2$ are the results of Liebeck and Seitz \cite{LS} which classify the nilpotent classes
and also compute their centralizers in the algebraic group.   The descriptions given in even characteristic require some interesting 
combinatorial analysis to get  nice generating functions and recognize them  as related to partitions.   See also \cite{Xu} which also counts the number of nilpotent orbits in characteristic $2$ for the orthogonal and symplectic Lie algebras over a finite field.  The methods there are very different and use intersection cohomology and local systems.  

 Let $p(n)$ denote the partition function, that is the number of positive integer partitions of $n$.   Let $P(u) = \sum p(n)u^n$. 
 
 We first state our results for an algebraically closed field $k$.     For $\mathfrak{gl}(n,k)$, the number of nilpotent classes is $p(n)$.  
 
 The next result is Theorem \ref{countnil} and Lemma \ref{countnil2}.

 \begin{thm}  Let $k$ be an algebraically closed field of characteristic $2$.  
 \begin{enumerate}
 \item  The generating function for the number of nilpotent orbits in $\mathfrak{sp}(2n,k)$ is $P(u)^2$. 
 \item  The generating function for the number of orthogonal group nilpotent orbits in $\mathfrak{so}(2n,k)$ is $(1-u)P(u)^2$.
 \end{enumerate} 
 \end{thm} 
 
 In the finite case, the answer is the same as for the algebraically closed case for $\mathfrak{sp}$ but is not for $\mathfrak{so}$.  
 See  Theorem \ref{soeven}. 
 
\begin{thm} \label{t:nilpotent}   Let $q$ be a power of $2$. 
\begin{enumerate} 
\item   The number of nilpotent orbits in $\mathfrak{sp}(2n,q)$ is $\sum_{i=0}^n p(i)p(n-i)$. 
\item   If $n$ is odd, the number of nilpotent orbits under the action of  $O^{\epsilon}(2n,q)$  in $\mathfrak{so}^{\epsilon}(2n,q)$ is $(1/2)\sum_{i=0}^n p(i)p(n-i)$.
\item if $n$ is even, the number of nilpotent orbits under the action of $O^{\epsilon}(2n,q)$ in $\mathfrak{so}^{\epsilon}(2n,q)$ is $(1/2) (\sum_{i=0}^n p(i)p(n-i) + \epsilon \cdot p(n/2))$.
\end{enumerate}
\end{thm}

We have a similar result for the total number of classes for $\mathfrak{sp}$ and $\mathfrak{so}$.  See Theorems \ref{finegen},  \ref{sum} and \ref{diff}.  In particular, 
we see, letting $n$ denote the rank, that the total number of orbits in the sympectic case is close to $(n+1)q^n$ and in the orthogonal case is close to $(n+1)q^n/2$.

For semisimple classes, the most interesting case is that of the orthogonal Lie algebra.  We find that the number of semisimple classes is much larger than
in the other cases.  Moreover, there is a very nice exact formula.  

\begin{thm} \label{ss}  Let $q$ be a power of $2$. For $n \geq 1$, the number of orbits of $O^{\epsilon}(2n,q)$ on semisimple elements in $\mathfrak{so}^{\epsilon}(2n,q)$ is 
$((n+1)q^n   -  (n-1)q^{n-1})/2$.  
\end{thm} 

Note that by Burnside's Lemma (due to Frobenius),  we see that if $G$ is a finite group of Lie type acting on its Lie algebra $\mathfrak{g}$,  then the proportion of commuting pairs 
of elements in $G \times \mathfrak{g}$ is $k(\mathfrak{g})/|\mathfrak{g}|$, where $k(\mathfrak{g})$ is the number of $G$ orbits on $\mathfrak{g}$.  Thus our estimates on $k(\mathfrak{g})$ give estimates for this commuting probability.  We see that this is 
related to properties of the commuting subvariety of $X \times \mathfrak{x}$ where $X$ is the corresponding algebraic group and $\mathfrak{x}$ is the Lie algebra of $X$.  We discuss this in detail in the final section.

In the next three sections, we consider the various cases (linear and unitary groups, symplectic groups and orthogonal groups).  

Diaconis brought to our attention the paper \cite{DZ}, which combines Markov chain theory (the Burnside process) with importance sampling to estimate the number of conjugacy classes of the unitriangular group over a finite field, giving a new approach to a conjecture of Higman. This is not directly useful for us as we have exact generating functions, but their approach does work to some extent for other group actions, and means that one potential use of our results (and also the results of \cite{FG}) is to serve as a testing ground for the approach of \cite{DZ}. The companion paper \cite{DH} uses the Burnside process to rapidly generate an approximately uniform conjugacy class of the symmetric group (a random partition) and could potentially be extended to our setting.
  
\section{$GL(n,q)$ and $GU(n,q)$}

Recall that the Lie algebra $\mathfrak{gl}(n,q)$ of $GL(n,q)$ consists of all $n \times n$ matrices over $F_q$.

The following result is Proposition 1.10.2 of Stanley \cite{St}.

\begin{thm} \label{fromstan} The number of orbits of $GL(n,q)$ in the adjoint action on its Lie algebra is equal to the coefficient
of $u^n$ in \begin{equation} \label{firstgen} \prod_{i \geq 1} \frac{1}{1-qu^i}. \end{equation}
\end{thm}

{\it Remark:} Setting $q=1$ in \eqref{firstgen} gives the generating function for the number of conjugacy classes in the symmetric group,
which is the Weyl group of $GL(n,q)$.  As we will see later, the analogous result does not hold for odd characteristic symplectic
groups (though it does hold for even characteristic symplectic groups). However setting $q=1$ in \eqref{firstgen} also gives the generating function for the number of nilpotent orbits of $GL(n,q)$
on its Lie algebra, and the analogous result does hold for the symplectic groups (in both odd and even characteristic).

\[ \]
A result of Maslen and Rockmore \cite{MR} is that the number of conjugacy classes of $GL(n,q)$ is at most $q^n$. The number of conjugacy
classes of $GL(n,q)$ is clearly less than the number of orbits of $GL(n,q)$ on its Lie algebra, but as the next result shows, the number of orbits on
the Lie algebra is still bounded by $C q^n$, where $C$ is a universal constant.  It also shows that as $q \rightarrow \infty$ 
(independently of $n$) that the number of classes divided by $q^n$ approaches $1$.  On the other hand, for $n \geq 2$, the number is strictly greater 
than $q^n$ since there are already $q^n$ distinct characteristic polynomials.  

\begin{thm} \label{glbound}
\begin{enumerate}
\item The number of orbits of $GL(n,q)$ on its Lie algebra is at most $3.5 q^n$.
\item The number of orbits of $GL(n,q)$ on its Lie algebra is at most $q^n + 5q^{n-1}$.
\end{enumerate}
\end{thm}

\begin{proof} For the first assertion, write \eqref{firstgen} as
\[ \prod_{i \geq 1} \frac{1-u^i}{1-qu^i} \prod_{i \geq 1} \frac{1}{1-u^i}. \]  It follows from Theorem \ref{fromstan}
that the number of orbits of $GL(n,q)$ on its Lie algebra is equal to
\[ \sum_{m=0}^n \left( Coef. \ u^{n-m} \ in \prod_{i \geq 1} \frac{1-u^i}{1-qu^i} \right)
\left( Coef. \ u^m \ in \prod_{i \geq 1} \frac{1}{1-u^i} \right).\] 

Now \[ \prod_{i \geq 1} \frac{1-u^i}{1-qu^i} \] is the generating function for the number of conjugacy classes of $GL(n,q)$,
so by the result of Maslen and Rockmore,
\[ Coef. \ u^{n-m} \ in  \prod_{i \geq 1} \frac{1-u^i}{1-qu^i} \leq q^{n-m}.\]

It follows that the number of orbits of $GL(n,q)$ on its Lie algebra is at most
\begin{eqnarray*}
\sum_{m=0}^n q^{n-m} Coef. \ u^m \ in \ \prod_{i \geq 1} \frac{1}{1-u^i}
& \leq & \sum_{m=0}^{\infty} q^{n-m} Coef. \ u^m \ in \ \prod_{i \geq 1} \frac{1}{1-u^i} \\
& = & q^n \prod_{i \geq 1} \frac{1}{1-1/q^i} \\
& \leq & 3.5 q^n,
\end{eqnarray*} proving the first assertion. The second assertion follows since
\[  \prod_{i \geq 1} \frac{1}{1-1/q^i} \leq 1+5/q.\]
\end{proof}

{\it Remark:} The bounds on the infinite product $\prod_{i \geq 1} (1-1/q^i)$ used in the previous theorem are immediate
from Euler's pentagonal number theorem:

\[ \prod_{i \geq 1} (1-x^i) = 1-x-x^2+x^5+x^7-x^{12}-x^{15}+x^{22}+x^{26} - \cdots.\]

\[ \]

{\it Remark:} One can also derive the generating function for the number of the  orbits on $\mathfrak{gl}(n,q)$ by noting that it is the product of the
generating function for the number of classes of $GL(n,q)$ and the generating function for the partition function.  

We now deal with the Lie algebra $\mathfrak{gu}(n,q)$ of $GU(n,q)$.  We show:

\begin{thm} \label{t:gu}  The number of orbits of $GU(n,q)$ on its Lie algebra is the same as the number of orbits of $GL(n,q)$ on its Lie algebra.
\end{thm}

\begin{proof}  Note that the Lie algebra  of $GU(n,q)$ is the space of skew Hermitian $n \times n$ matrices in the space of $n \times n$ matrices
over the field of $q^2$ elements.  To count orbits, we can work with Hermitian matrices since the number of orbits are the same (one space is a scalar 
multiple of the other).   

Suppose that $A$ is Hermitian.  Then $\sigma(A) = A^{\top}$ where $\sigma$ is the $q$th power automorphism.  Since $A$ and $A^{\top}$ are similar,
it follows that $\sigma(A)$ and $A$ are similar and so all the invariant factors of $A$ are defined over the field of $q$ elements and so $A$ is similar to
a matrix $B \in \mathfrak{gl}(n,q)$.  Conversely, if $B \in \mathfrak{gl}(n,q)$, then $B$ is self adjoint with respect to some invertible symmetric bilinear form over 
the field of $q$-elements  (cf. \cite{TZ}) and so also with respect to an invertible Hermitian form over the field of $q^2$ elements. In particular, 
 $B$ is conjugate to a Hermitian matrix.   It only remains to show that if $A$ and $B$ are Hermitian matrices which are conjugate over the algebraic closure, then
are conjugate by an element of $GU(n,q)$.   This follows by the Lang-Steinberg theorem using the fact that centralizers of matrices in $GL(n,k)$ are connected
(the centralizer in $GL(n,k)$ is the group of units of the algebra centralizer; here $k$ is an algebraically closed field containing the field of $q$ elements).
\end{proof}  

Thus, the results for the action of $GL(n,q)$ acting on $\mathfrak{gl}(n,q)$ apply and we get the same estimates (and the same generating function).  \\

{\it Remark:}   If we want to deal with $\mathfrak{sl}(n,q)$ and $\mathfrak{su}(n,q)$, the results are similar but slightly different and the answers for $\mathfrak{sl}(n,q)$ and
$\mathfrak{su}(n,q)$ are not the same.  

\section{Symplectic groups}

This section treats the finite symplectic groups. The behavior is very different in odd and even characteristic, so we treat these separately.

\subsection{Odd characteristic} 

Recall that $\mathfrak{sp}(2n,q)$ is the space of skew adjoint operators on the natural $2n$-dimensional space $V$ with respect to the alternating form
defining the symplectic group.  

The orbit of an element is determined by a bit more than its rational canonical form.   The semisimple part of the element is determined up to conjugacy by its characteristic polynomial.  

In what follows, for $p(x)$ a monic degree $d$ polynomial, we define $p^*(x)=(-1)^d p(-x)$. We call $p(x)$ self-negative if $p(x)=p(-x)$ and self-conjugate if $p^*(x)=p(x)$. 

We note the following  well known elementary linear algebra result.  This follows easily from the fact that for a semisimple skew adjoint operator,  the $0$ eigenspace  is nondegenerate and 
all other eigenspaces are totally singular, and the $\alpha$ and $\beta$ eigenspaces are orthogonal unless $\beta = - \alpha$. 

\begin{lem} \label{l:ss-sp}   Let $T$ be a skew adjoint operator on $V$ with respect to the alternating form.  Then the characteristic polynomial of $T$
is a product of factors which are of the form $x^2$, $f(x)$ where $f(x)$ is a monic irreducible polynomial with $f(x)=f(-x)$, and $g(x)g^*(x)$ where $g(x)$ is a monic irreducible polynomial such that $g(x) \neq g^*(x)$.
\end{lem}  

The following elementary lemma will be helpful.

\begin{lem} \label{easy}
\begin{enumerate}
\item If $deg(p(x))$ is even, then $p(x)$ is self-negative if and only if it is self-conjugate.
\item If $p(x)$ is irreducible and self-conjugate, then it has even degree or $p(x)=x$.
\end{enumerate}
\end{lem}

\begin{proof} The first assertion is obvious. For the second assertion, if $p(x)$ is self-conjugate then if $a$ is a root, so is $-a$. So if $deg(p(x))$ is odd, it must have $0$ as a root, so $p(x)=x$.
\end{proof}

An orbit of $Sp(2n,q)$ on its Lie algebra is determined by  its canonical form and a bit more.

We split the orbit up by the irreducible factors of its characteristic polynomial. There are three cases as noted in Lemma \ref{l:ss-sp}.  
So let $A$ be a semisimple element with characteristic polynomial $f(x)$.  \\ 

{\bf Case 1}: $f(x)=x^{m}$.  The centralizer of $A$ in this case is just $\mathfrak{sp}(m,q)$ and  so the number of orbits on elements with given semisimple part  is just the number of nilpotent classes in $\mathfrak{sp}(m,q)$.  In this case, it is well known that the  orbits correspond to partitions of $m$ where the odd parts occur with even multiplicity, and for each even 
size that occurs there is a choice of sign.  See \cite[Chapter 3]{LS}.  \\

{\bf Case 2}: $f(x) = g(x)^m$ with $g(x)$ monic irreducible and self-negative. Then the degree $d$ of $g$ is even. The centralizer of $A$ is $\mathfrak{gu}(m,q^{d/2})$. and the number of nilpotent classes commuting with the semisimple element is the number of partitions of $m$. \\

{\bf Case 3}: $f(x)=(g(x) g^*(x))^m$, where $g(x)$ is monic irreducible of degree $d$ and $g(x) \neq g^*(x)$. Suppose that $f(x)$ is the largest factor dividing the characteristic polynomial. Then the centralizer of the semisimple part is $\mathfrak{gl}(m, q^d)$.  Thus the nilpotent classes centralizing the semisimple element correspond to partitions of  $m$.  \\

We can now put these together to obtain a generating function for the number of orbits. 

In what follows, let $N^*(q;d)$ denote the number of monic irreducible self-conjugate polynomials of degree $d$, and let $M^*(q;d)$ denote the number of unordered pairs $\{g(x),g^*(x)\}$ such that  $g(x)$ is monic, irreducible of degree $d$ and $g(x) \neq g^*(x)$.

The following lemma will be helpful.

\begin{lem} \label{counterrep}
\[ \prod_{d \geq 1} (1-u^d)^{-N^*(q;2d)} \prod_{d \geq 1} (1-u^d)^{-M^*(q;d)} = \frac{1-u}{1-qu}.\]
\end{lem}

\begin{proof} Consider the coefficient of $u^{2n}$ in 
\begin{equation} \label{infprod} \prod_{d \geq 1} (1-u^{2d})^{-N^*(q;2d)} \prod_{d \geq 1} (1-u^{2d})^{-M^*(q;d)} .
\end{equation}

Now by part 2 of Lemma \ref{easy}, if $p(x)$ is monic self-conjugate and has non-zero constant term, it factors uniquely into a product of even degree monic irreducible self-conjugate polynomials $f_i(x)$ and pairs $g_i(x) g_i^*(x)$ where $g_i(x)$ is monic irreducible and non-self conjugate. So the coefficient of $u^{2n}$ in \eqref{infprod} is equal to the total number of monic self-conjugate polynomials of degree $2n$ with non-zero constant term. By part 1 of Lemma \ref{easy}, this is the number of monic self-negative polynomials of degree $2n$ with non-zero constant term, which is $1$ if $n=0$ and $q^n-q^{n-1}$ if $n \geq 1$. Hence \eqref{infprod} is equal to $\frac{1-u^2}{1-u^2q}$ and the lemma follows. \end{proof}

First we record a generating function for the number of nilpotent orbits of $Sp(2n,q)$ on its Lie algebra.

\begin{thm} \label{nilcount} Suppose that the characteristic is odd. Then the number of nilpotent orbits of $Sp(2n,q)$ on its Lie algebra is equal to the coefficient of $u^n$ in \[ A(u):= \prod_{i \geq 1} \frac{(1+u^i)^2}{1-u^i}.\]
\end{thm}

\begin{proof} By Case 1 above, the number of nilpotent orbits is equal to the coefficient of $u^{2n}$ in
\begin{eqnarray*}
& & \prod_{i \ odd} \frac{1}{1-u^{2i}} \prod_{i \ even} (1+2u^i+2u^{2i}+2u^{3i}+\cdots) \\
& = & \prod_{i \ odd} \frac{1}{1-u^{2i}} \prod_{i \ even} \frac{1+u^i}{1-u^i}\\
& = & \prod_{i \ odd} \frac{1}{1-u^{2i}} \prod_{i \geq 1} \frac{1+u^{2i}}{1-u^{2i}} \\
& = & \prod_{i \ odd} \frac{1}{1-u^{2i}} \prod_{i \geq 1} \frac{(1+u^{2i})^2}{1-u^{4i}}\\
& = & \prod_{i \geq 1} \frac{(1+u^{2i})^2}{1-u^{2i}}.
\end{eqnarray*} \end{proof}

Theorem \ref{sympoddgen} derives a generating function for the number of orbits of $Sp(2n,q)$ on its Lie algebra.

\begin{thm} \label{sympoddgen} Suppose that the characteristic is odd. Then the number of orbits of $Sp(2n,q)$ on its Lie algebra is
equal to the coefficient of $u^n$ in
\begin{equation} \label{usegen} 
\prod_{i \geq 1} \frac{(1+u^i)^2}{1-qu^i}.
\end{equation}
\end{thm} 

\begin{proof} By the above discussion and the previous theorem, the total number of orbits is the coefficient of $u^n$ in
\begin{equation} \label{t1} A(u) \prod_{d \geq 1} (\sum_{\lambda} u^{d |\lambda|})^{N^*(q;2d)} \prod_{d \geq 1}  (\sum_{\lambda} u^{d |\lambda|})^{M^*(q;d)}, \end{equation} where the sums are over all integer partitions $\lambda$.

This is equal to \[ A(u) \prod_{d \geq 1} \prod_{i \geq 1} (1-u^{id})^{-N^*(q;2d)} \prod_{d \geq 1} \prod_{i \geq 1} (1-u^{id})^{-M^*(q;d)}.\]
Switching the order of the infinite products and applying Lemma \ref{counterrep} yields
\[ A(u) \prod_{i \geq 1} \prod_{d \geq 1}  (1-u^{id})^{-N^*(q;2d)} (1-u^{id})^{-M^*(q;d)} = A(u) \prod_{i \geq 1} \frac{1-u^i}{1-qu^i}.\] This completes the proof. \end{proof}

{\it Remark:} Setting $q=1$ in $(\ref{usegen})$ gives the generating function for the number of nilpotent orbits. 

\[ \]
Next we give an upper bound for the number of orbits of $Sp(2n,q)$ on its Lie algebra.  Note that a trivial lower bound of $q^n$ comes by considering semisimple classes
(the number of possible characteristic polynomials).  

\begin{thm} Suppose that the characteristic is odd.
\begin{enumerate}
\item The number of orbits of $Sp(2n,q)$ on its Lie algebra is at most $4.5q^n$.
\item The number of orbits of $Sp(2n,q)$ on its Lie algebra is at most $q^n + 11q^{n-1}$.
\end{enumerate}
\end{thm}

\begin{proof} For the first part, we argue as in Theorem \ref{glbound}. From \eqref{usegen}, the number of orbits is equal to
\[ \sum_{m=0}^n \left( Coef. \ u^{n-m} \ in \ \prod_i \frac{1-u^i}{1-qu^i} \right) \left( Coef. u^m \ in \ \prod_i \frac{(1+u^i)^2}{(1-u^i)} \right).\]
As in the proof of Theorem \ref{glbound}, the coefficient of $u^{n-m}$ in
\[ \prod_i \frac{1-u^i}{1-qu^i} \] counts the number of conjugacy classes of $GL(n-m,q)$ and is at most $q^{n-m}$.

It follows that the number of orbits of $Sp(2n,q)$ on its Lie algebra is at most
\begin{eqnarray*}
q^n \sum_{m=0}^n \frac{1}{q^m} Coef. \ u^m \ in \ \prod_i \frac{(1+u^i)^2}{(1-u^i)}
& \leq &  q^n \sum_{m=0}^{\infty} \frac{1}{q^m} Coef. \ u^m \ in \ \prod_i \frac{(1+u^i)^2}{(1-u^i)} \\
& = & q^n \prod_{i \geq 1} \frac{(1+1/q^i)^2}{(1-1/q^i)}\\
& \leq & 4.5 q^n
\end{eqnarray*}

The second assertion follows from the proof of the first assertion and the inequality
\[ \prod_{i \geq 1} \frac{(1+1/q^i)^2}{(1-1/q^i)} \leq 1 + \frac{11}{q}.\]
\end{proof}

\subsection{Even characteristic}

To begin we enumerate the number of nilpotent orbits of $Sp(2n,q)$ on its Lie algebra.    The following result says that the number of nilpotent 
orbits in $\mathfrak{sp}(2n,q)$ is $\sum_{i+j=n} p(i)p(j)$ where $p(j)$  is the partition function (equivalently the generating function for
the number of the orbits is the square of the generating function for the partition function).   This is equal to the number of conjugacy classes in the Weyl
group of type $B_n$.   It follows by \cite[5.12]{LS} that the stabilizers  of nilpotent elements in $\mathfrak{sp}(2n,k)$ with $k$ an algebraically closed 
field of even characteristic are connected and so all classes are defined over the prime field and remain a single class over the finite field. 

The description of the classes given in \cite{LS} is fairly complicated and it seems surprising that the answer comes out so nicely.   We also note that
by \cite{LS}, it follows that the nilpotent orbits are all defined over $F_2$ and moreover the stabilizers in the algebraic group are connected.  Thus,
the number of orbits is the same whether counting over a finite field or the algebraic closure by the Lang-Steinberg theorem.

\begin{thm} \label{countnil} Suppose that the characteristic is even. Then the number of nilpotent orbits of $Sp(2n,q)$ on its Lie algebra is equal to
the coefficient of $u^n$ in \[ \prod_{i \geq 1} \frac{1}{(1-u^i)^2}.\]
\end{thm}

{\it Remark:} Let us say a bit about the asymptotics of $a_n$, the number of nilpotent orbits of $Sp(2n,q)$ on its Lie algebra.
From the Hardy-Littlewood-Rademacher asymptotics for the
partition function $p_n$, it is a routine computation to show that
$$a_n=\sum_{k=0}^n p_k p_{n-k} \sim \frac{1}{4\cdot 3^{3/4}\cdot n^{5/4}} e^{2\pi \sqrt{n/3}} \left(1 - \frac{4\pi^2+135}{48\pi \sqrt{3n}} + O\left(\frac{1}{n}\right)\right).$$
One could carry this further to get a series in powers of $1/\sqrt n$, and if
you look at the OEIS page for this sequence it gives the $1/n$ term
explicitly. The Hardy-Littlewood-Rademacher formula contains lower order terms 
(which are so small that they are negligible compared to every term in the
asymptotic series above) with
periodic behavior (with every possible period), so turning the asymptotic
formula into a hard inequality is a little tricky. Based on numerical evidence
only, one gets that if one writes
$$a_n=\frac{1}{4\cdot 3^{3/4}\cdot n^{5/4}} e^{2\pi \sqrt{n/3}} \left(1 - \frac{\delta_n}{\sqrt{n}}\right),$$
then for all $n\ge 1$
$$\frac{1}{2} < \delta_n < \frac{4\pi^2+135}{48\pi \sqrt{3}}=0.668\dots.$$
If one is committed, proving this should be routine, but very tedious, and
weakening the constants it would give the hard bounds
$$\frac{1}{4\cdot 3^{3/4}\cdot n^{5/4}} e^{2\pi \sqrt{n/3}} \left(1 - \frac{1}{\sqrt{n}}\right) < a_n < \frac{1}{4\cdot 3^{3/4}\cdot n^{5/4}} e^{2\pi \sqrt{n/3}}$$
for all $n\ge 1$. From this it follows that
$$d_n = a_n - a_{n-1} \sim \frac{\pi}{12\cdot 3^{1/4}\cdot n^{7/4}} e^{2\pi \sqrt{n/3}} \left(1-\frac{7(4\pi^2+45)}{48\pi\sqrt {3n}}+O\left(\frac{1}{n}\right)\right),$$
and based only on numeric evidence (and using weak constants)
$$\frac{\pi}{12\cdot 3^{1/4}\cdot n^{7/4}} e^{2\pi \sqrt{n/3}} \left(1-\frac{3}{\sqrt n}\right) < d_n < \frac{\pi}{12\cdot 3^{1/4}\cdot n^{7/4}} e^{2\pi \sqrt{n/3}}$$
for all $n\ge 1$.

\begin{proof} (Of Theorem \ref{countnil})
We use the results and notation of Liebeck and Seitz. Let $k$ be the algebraic closure of the field of $2$ elements.   According to \cite[Lemma 5.4]{LS},  every nilpotent element in $\mathfrak{sp}(2n,k)$
is conjugate to a unique element that is of the form $\sum_i W(m_i)^{a_i} + \sum_i W_{\ell_i}(n_i) + \sum_j V(k_j)$ subject to various constraints. For definitions of the indecomposables $V(m)$, $W(m)$, and $W_l(m)$ see \cite[Section 5.1]{LS}. A combinatorial description of the constraints will be given below. In particular, $\sum_i a_im_i + \sum n_i + \sum k_i = n$ (each building block corresponds to matrix of twice the size of the parameter).  Moreover, \cite[Theorem 5.12]{LS}, the centralizer of any nilpotent element of $\mathfrak{sp}(2n,k)$ is connected and thus we have precisely the same description of nilpotent classes over any finite field of characteristic $2$.

We can convert this to a purely combinatorial problem as follows. 
 
We consider three classes of building blocks. We have

\begin{enumerate}
\item For each positive integer $m$, we have an object $W(m)$ of weight $m$. We
can use any number of copies of these objects.

\item For each pair of positive integers $(\ell,n)$ with $0<\ell<n/2$, we have
an object $W_{\ell}(n)$ of weight $n$. We can use multiple objects of this
type, but the $n$-values corresponding to them must be distinct. Furthermore, 
if we use objects with indices $(\ell_i,n_i)$, where we arrange them so that
the sequence $(n_i)$ is in (strictly) decreasing order, then the sequences
$(\ell_i)$ and $(n_i-\ell_i)$ are also strictly decreasing.

\item For each positive integer $k$, we have an object $V(k)$ of weight $k$. We
can use at most two of each of these objects.
\end{enumerate}

\noindent
Finally there is an interaction between objects of types (2) and (3). If the
sequence of parameters for the objects of types (2) and (3) used are
$(n_i)$, $(\ell_i)$, and $(k_j)$, respectively, then

\bigskip

(4) For each pair $i,j$, we have either $k_j>n_i-\ell_i$ or $k_j<\ell_i$.

\bigskip

Let $a_n$ be the number of ways to choose a collection of objects subject to
these conditions with total weight $n$. Let $b_n$ be the same number where
we only allow the last two classes of objects, and let
$$A(q)=\sum_{n=0}^{\infty} a_n q^n,\; {\rm and}\; B(q)=\sum_{n=0}^{\infty} b_nq^n$$
be their ordinary generating functions. Since there is no interaction between
objects of type (1) and the other two types, and objects of type (1) allow us
to take an arbitrary sum of positive integers, hence a partition, we have
$$A(q) = P(q) B(q),$$
where $P(q) = \prod_{n=1}^{\infty} \frac{1}{1-q^n}$ is the partition function.
We prove that in fact $B(q)=P(q)$ and hence
$$A(q)=P(q)^2 = 1 + 2 q + 5 q^2 + 10 q^3 + 20 q^4 + 36 q^5 + 65 q^6 + \dots.$$
This makes $(a_n)$ the sequence A000712 from the OEIS (online encyclopedia of
integer sequences).

Before we begin the real proof, let us quickly discuss the other two pieces.
Pieces of type (3) on their own are relatively straightforward. They are
partitions where we are restricted to at most 2 copies of each part. The
generating function for the number of them is just
$$\prod_{k=1}^{\infty} (1+q^k+q^{2k}) = \prod_{k=1}^{\infty} \frac{1-q^{3k}}{1-q^k} = \frac{P(q)}{P(q^3)}.$$

The objects of type (2) are a little more subtle.
If we define the sequence $\ell'_i=n_i-\ell_i$, then we have $0<\ell_i<\ell'_i$
for all $i$ and the sequences $(\ell_i)$ and $(\ell'_i)$ are strictly
decreasing. Thus putting the $\ell'_i$ on the top row and the $\ell_i$ on the
bottom row, we can regard these as giving a way to fill a $2\times r$
grid with positive integers (for some $r$), so that they are strictly
decreasing as we go either down or right.
Let $G_r$ be the set of all such $2\times r$ grids. Conversely, if we take
such a grid $g\in G_r$ and use its rows to define $\ell_i$ and $\ell'_i$, then
$n_i=\ell_i+\ell'_i$ will be strictly decreasing, we will have
$2\ell_i<\ell_i+\ell'_i=n_i$, and the total weight $\sum_i n_i$ will be just the
sum $|g|$ of all the entries of $g$. Though it is irrelevant to the proof 
below, these grids are up to cosmetic changes
plane partitions based on a $2\times r$ grid, and the generating function for
the number of them with a fixed sum of entries is
$$\sum_{g\in G_r} q^{|g|}=\frac{q^{r(r+2)}}{(1-q)(1-q^2)^2\cdots (1-q^r)^2(1-q^{r+1})}.$$
We can view the columns of $g$ as defining
closed intervals $[\ell_i,\ell'_i]$ and define the excluded set for $g$ to be
$$X_g = \bigcup_i\; [\ell_i,\ell'_i].$$
Then the interaction condition between objects of types (2) and (3) can be
stated as saying that if we fix $g$, then $k_j\notin X_g$. Hence we should omit
the corresponding factors $\prod_{k\in X_g} (1+q^k+q^{2k})$ from the generating
function above. Alternately, we can agree to cancel them off, which leads to
looking at the generating function
$$C(q) = 1 + \sum_{r=1}^{\infty} \sum_{g\in G_r} \frac{q^{|g|}}{\prod_{k\in X_g} (1+q^k+q^{2k})},$$
and with this definition we have
$$B(q) = \frac{P(q)}{P(q^3)} \cdot C(q),$$
and hence our goal is to show that $C(q)=P(q^3)$.

\bigskip

Now the proof is just sorting out this sum. Though this sum is a little
elaborate, we can break it down inductively. For $k<n$, let $C_{k,n}$ be the
sum of all terms from $C(q)$ where the first column of the grid is 
$\left(\begin{matrix} n\cr k\end{matrix}\right)$. So for example, if $k=1$,
then the only such grid is the $2\times 1$ grid with just this column and hence
$$C_{1,n}=\frac{q^{n+1}}{\prod_{i=1}^n (1+q^i+q^{2i})}.$$
With this as our base case, we will prove by induction on $k$ that
$$C_{k,n}=\frac{q^{n+k} \prod_{i=n+1}^{n+k-1} (1-q^i)}{\prod_{i=1}^{k-1} (1-q^i)\cdot \prod_{i=1}^n (1+q^i+q^{2i})}.$$
If we have a grid with this first column (and it is not just the $2\times 1$
grid with only this column), and we peel off the first column, then we will be
left with a smaller grid $g'$, say with first column corresponding to $(j,m)$.
We must have $j<k$ and $j<m<n$ (and conversely any such grid $g'$ arises), and
hence
$$C_{k,n}=\frac{q^{n+k}}{\prod_{i=k}^n (1+q^i+q^{2i})} + \sum_{j=1}^{k-1} \sum_{m=j+1}^{n-1} \frac{q^{n+k} C_{j,m}}{\prod_{i=\max(k,m+1)}^n (1+q^i+q^{2i})},$$
where we have used the obvious fact that $|g|=n+k+|g'|$ and we have noted
that the interval $[\max(k,m+1),n]=[k,n]\backslash [j,m]$ is always exactly
$X_g\backslash X_{g'}$.
Since we have already seen the base case, all that remains is to check that 
our proposed formula satisfies this recursion. Hence we only need to prove
$$\frac{q^{n+k} \prod_{i=n+1}^{n+k-1} (1-q^i)}{\prod_{i=1}^{k-1} (1-q^i)\cdot \prod_{i=1}^n (1+q^i+q^{2i})}
=\frac{q^{n+k}}{\prod_{i=k}^n (1+q^i+q^{2i})}$$
$$+ \sum_{j=1}^{k-1} \sum_{m=j+1}^{n-1} \frac{q^{n+k+m+j} \prod_{i=m+1}^{m+j-1} (1-q^i)}{\prod_{i=1}^{j-1} (1-q^i)\cdot \prod_{i=1}^m (1+q^i+q^{2i}) \prod_{i=\max(k,m+1)}^n (1+q^i+q^{2i})}.$$
After clearing denominators and cancelling a common factor, this means we need
to prove
$$\prod_{i=n+1}^{n+k-1} (1-q^i) = \prod_{i=1}^{k-1} (1-q^{3i})
+ \sum_{j=1}^{k-1} \sum_{m=j+1}^{n-1} q^{m+j} \prod_{i=m+1}^{m+j-1} (1-q^i)\cdot \prod_{i=j}^{k-1} (1-q^i)\cdot \prod_{i=m+1}^{k-1} (1+q^i+q^{2i}),$$
where the final product on the right should be interpreted as an empty product,
hence equal to $1$, if $m+1>k-1$.

To prove this it suffices to show that it holds for $n=k+1$ and that the 
identity one gets by taking a difference of both sides with respect to $n$
holds. This difference formula reads
$$\prod_{i=n+2}^{n+k} (1-q^i) - \prod_{i=n+1}^{n+k-1} (1-q^i) =
\sum_{j=1}^{k-1} q^{n+j} \prod_{i=n+1}^{n+j-1} (1-q^i)\cdot \prod_{i=j}^{k-1} (1-q^i),$$
which, since the left-hand side factors as $q^{n+1}(1-q^{k-1})\prod_{i=n+2}^{n+k-1} (1-q^i)$, cancels down to 
$$\prod_{i=n+2}^{n+k-1} (1-q^i) =
\sum_{j=1}^{k-1} q^{j-1} \prod_{i=n+1}^{n+j-1} (1-q^i)\cdot \prod_{i=j}^{k-2} (1-q^i).$$
Writing
$$q^{j-1} (1-q^{n+1})=(1-q^{n+j}) - (1-q^{j-1}),$$
we see that the summand on the right is
$$\prod_{i=n+2}^{n+j} (1-q^i)\cdot \prod_{i=j}^{k-2} (1-q^i) - \prod_{i=n+2}^{n+j-1} (1-q^i)\cdot \prod_{i=j-1}^{k-2} (1-q^i),$$
so this sum telescopes giving exactly this.

Thus we only need to check the case $n=k+1$, which reads
$$\prod_{i=k+2}^{2k} (1-q^i) = \prod_{i=1}^{k-1} (1-q^{3i})
+ \sum_{j=1}^{k-1} \sum_{m=j+1}^k q^{m+j} \prod_{i=m+1}^{m+j-1} (1-q^i)\cdot \prod_{i=j}^{k-1} (1-q^i)\cdot \prod_{i=m+1}^{k-1} (1+q^i+q^{2i}).$$
Pulling out the $m=k$ term (which is weird because of the empty product), and
regrouping some products, this becomes
\begin{equation*}\begin{split}
\prod_{i=k+2}^{2k} (1-q^i) = \prod_{i=1}^{k-1} (1-q^{3i})
& +  \sum_{j=1}^{k-1} \frac{q^{k+j}}{1-q^k} \prod_{i=j}^{k+j-1} (1-q^i)\\
& + \sum_{m=2}^{k-1} \sum_{j=1}^{m-1} q^{m+j} \prod_{i=j}^{m+j-1} (1-q^i)\cdot \prod_{i=m+1}^{k-1} (1-q^{3i}).
\end{split}\end{equation*}
Since
$$q^{m+j} \prod_{i=j}^{m+j-1} (1-q^i)=\frac{q^{m+1}}{1-q^{m+1}}\left[\prod_{i=j}^{m+j} (1-q^i)-\prod_{i=j-1}^{m+j-1} (1-q^i)\right],$$
(and using this with $m\mapsto k$ in the first sum), we see that the $j$ sums
telescope, so we reduce to
\begin{equation*}\begin{split}
\prod_{i=k+2}^{2k} (1-q^i) = \prod_{i=1}^{k-1} (1-q^{3i})
& + \frac{q^{k+1}}{(1-q^k)(1-q^{k+1})} \prod_{i=k-1}^{2k-1} (1-q^i)\\
& + \sum_{m=2}^{k-1} \frac{q^{m+1}}{1-q^{m+1}} \prod_{i=m-1}^{2m-1} (1-q^i)\cdot \prod_{i=m+1}^{k-1} (1-q^{3i}),
\end{split}\end{equation*}
which simplifies to
$$\prod_{i=k+1}^{2k-1} (1-q^i) = \prod_{i=1}^{k-1} (1-q^{3i})
+ \sum_{m=2}^{k-1} \frac{q^{m+1}}{1-q^{m+1}} \prod_{i=m-1}^{2m-1} (1-q^i)\cdot \prod_{i=m+1}^{k-1} (1-q^{3i}),$$
This $m$ sum also telescopes since
\begin{equation*}\begin{split}
\prod_{i=m+2}^{2m+1}(1-q^i) - & (1-q^{3m}) \prod_{i=m+1}^{2m-1}(1-q^i)\\
&\quad =\left[(1-q^{2m})(1-q^{2m+1})-(1-q^{m+1})(1-q^{3m})\right] \prod_{i=m+2}^{2m-1}(1-q^i)\\
&\quad =q^{m+1}(1-q^{m-1})(1-q^m) \prod_{i=m+2}^{2m-1}(1-q^i),
\end{split}\end{equation*}
and this formula follows. Thus we have completed the inductive proof of our
asserted formula for $C_{k,n}$.

\bigskip

Plugging in this formula for $C_{k,n}$, our messy sum for $C(q)$ reduces to
just a double sum
$$C(q) = 1 + \sum_{n=2}^{\infty} \sum_{k=1}^{n-1} 
\frac{q^{n+k} \prod_{i=n+1}^{n+k-1} (1-q^i)}{\prod_{i=1}^{k-1} (1-q^i)\cdot \prod_{i=1}^n (1+q^i+q^{2i})}.$$
Using the identity
$$\frac{q^{k-1} \prod_{i=1}^{k-1} (1-aq^i)}{\prod_{i=1}^{k-1} (1-q^i)} =
\frac{\prod_{i=2}^{k} (1-aq^i)}{\prod_{i=1}^{k-1} (1-q^i)} 
- \frac{\prod_{i=2}^{k-1} (1-aq^i)}{\prod_{i=1}^{k-2} (1-q^i)},$$
with $a=q^n$, we see that the inner sum telescopes, giving
\begin{equation*}\begin{split}
C(q) & = 1 + \sum_{n=2}^{\infty} 
\frac{q^{n+1} \prod_{i=n+2}^{2n-1} (1-q^i)}{\prod_{i=1}^{n-2} (1-q^i)\cdot \prod_{i=1}^n (1+q^i+q^{2i})}\\
& =1+\sum_{n=2}^{\infty} \frac{q^{n+1} (1-q^{n-1})(1-q^n) \prod_{i=n+2}^{2n-1} (1-q^i)}{\prod_{i=1}^n (1-q^{3i})}.
\end{split}\end{equation*}
Using the same computation we just gave above (with $n$ in place of $m$)
$$q^{n+1}(1-q^{n-1})(1-q^n)= (1-q^{2n})(1-q^{2n+1})-(1-q^{n+1})(1-q^{3n}),$$
we see that this sum telescopes
\begin{equation*}\begin{split}
C(q) & = 1+\sum_{n=2}^{\infty}\left[
\frac{\prod_{i=n+2}^{2n+1} (1-q^i)}{\prod_{i=1}^n (1-q^{3i})}
- \frac{\prod_{i=n+1}^{2n-1} (1-q^i)}{\prod_{i=1}^{n-1} (1-q^{3i})}
\right]\\
& =\lim_{N\to \infty} \frac{\prod_{i=N+2}^{2N+1} (1-q^i)}{\prod_{i=1}^N (1-q^{3i})} = \prod_{i=1}^{\infty} \frac{1}{1-q^{3i}} = P(q^3).
\end{split}\end{equation*}.

\end{proof}

Armed with Theorem \ref{countnil} and the following lemmas, we will derive an extremely simple generating function for the number of orbits of $Sp(2n,q)$ on its Lie algebra in even characteristic.  Note that the characteristic polynomial of any element of $\mathfrak{sp}(2n,q)$ is a square and any monic degree $2n$ polynomial that is a square occurs as 
a characteristic polynomial. 

\begin{lem} \label{reduct} Let $f$ be a monic polynomial of degree $n$. Factor $f$ as a product of irreducibles, so
\[ f = \prod_{\phi} \phi^{m_{\phi}(f)} .\] Then the number of orbits of $Sp(2n,q)$ on its Lie algebra with characteristic polynomial $f^2$
is equal to the product over $\phi$ of the number of nilpotent orbits of $Sp(2m_{\phi}(f),q)$ on its Lie algebra. Moreover, all orbits
arise this way.
\end{lem}

\begin{proof}   Let $x$ be a semisimple element with characteristic polynomial $f^2$.  There is a unique class of such elements.  We need to count the number of 
nilpotent orbits in the centralizer of $x$.  Decomposing the space, we see it  suffices to consider the case that $f = \phi^m$ with $\phi$ irreducible.  
Then the centralizer of $x$ is $\mathfrak{sp}(2m,q^e)$ with $e$  the degree of $\phi$.  The number of nilpotent orbits is independent of $e$ and so the results holds.
\end{proof}

\begin{lem} \label{necklace} \[ \prod_{\phi} \frac{1}{1-u^{deg(\phi)}} = \frac{1}{1-uq}, \] where the product is over all monic irreducible polynomials with
coefficients in $F_q$ and $deg(\phi)$ denotes the degree of $\phi$.
\end{lem}

\begin{proof} By unique factorization, the coefficient of $u^n$ on the left hand side is equal to the total number of degree $n$ monic polynomials, which is
equal to $q^n$. The coefficient of $u^n$ on the right hand side is also equal to $q^n$.
\end{proof}

\begin{thm} \label{finegen} Suppose that the characteristic is even. Then the number of orbits of $Sp(2n,q)$ on its Lie algebra is equal to the coefficient of
$u^n$ in \begin{equation} \label{even} \prod_{i \geq 1} \left( \frac{1}{1-qu^i }\right)^2. \end{equation}
\end{thm}

\begin{proof} By Lemma \ref{reduct} and Theorem \ref{countnil}, the number of orbits is the coefficient of $u^n$ in
\[ \prod_{\phi} \prod_{i \geq 1} \left( \frac{1}{1-u^{i \cdot deg(\phi)}}  \right)^2, \]  where the product is over all monic irreducible
polynomials. Reversing the order of the products and applying Lemma \ref{necklace} gives that the number of orbits is the coefficient
of $u^n$ in \[ \prod_{i \geq 1} \left( \prod_{\phi} \frac{1}{1-u^{i \cdot deg(\phi)}} \right)^2 = \prod_{i \geq 1} \left( \frac{1}{1-qu^i} \right)^2, \]
as desired. \end{proof}

{\it Remark:} By Theorem \ref{countnil}, setting $q=1$ in \eqref{even} gives the number of nilpotent orbits of $Sp(2n,q)$ on its Lie algebra.
This is equal to the number of conjugacy classes of the Weyl group $B_n$.
\[ \]

Given our results for odd characteristic symplectic groups, and also the bounds on numbers of conjugacy classes of finite groups of Lie type \cite{FG},
one might expect that the number of orbits is at most a constant times $q^n$. Theorem \ref{goodbounds} shows this is not the case.

\begin{thm} \label{goodbounds} Suppose that the characteristic is even.
\begin{enumerate}
\item The number of orbits of $Sp(2n,q)$ on its Lie algebra is at least $(n+1)q^n$.
\item The number of orbits of $Sp(2n,q)$ on its Lie algebra is at most $12(n+1)q^n$.
\item The number of orbits of $Sp(2n,q)$ on its Lie algebra is at most \[(n+1) q^n (1+22/q).\]
\item For fixed $q$, the number of orbits of $Sp(2n,q)$ on its Lie algebra divided by $q^n$ is asymptotic to $C_q \cdot n$ as
$n \rightarrow \infty$ where \[ C_q = \prod_{i \geq 1} \left( \frac{1}{1-1/q^i} \right)^2.\]
\end{enumerate}
\end{thm}

\begin{proof} For part 1, taking the $i=1$ term in Theorem \ref{finegen} gives that the number of orbits is at least
$q^n$ multiplied by the coefficient of $u^n$ in $(1-u)^{-2}$, which is equal to $(n+1)q^n$.

For part 2, we write \eqref{even} as
\[ \prod_i \left( \frac{1-u^i}{1-qu^i}  \right)^2 \prod_i \left( \frac{1}{1-u^i} \right)^2,\] and argue as in the odd characteristic case.
Recalling that the coefficient of $u^n$ in \[ \prod_{i} \left( \frac{1-u^i}{1-qu^i} \right) \] is equal to the number of conjugacy classes
of $GL(n,q)$ which is at most $q^n$, it follows that the number of orbits is at most
\[ \sum_{m=0}^n \left( Coef. \ u^{n-m} \ in \ (1+qu+q^2u^2+\cdots)^2 \right) \left( Coef. \ u^m \ in \ \left( \prod_i \frac{1}{1-u^i} \right)^2 \right)\]
This is equal to
\[ \sum_{m=0}^n (n-m+1) q^{n-m}  Coef. \ u^m \ in \ \left( \prod_i \frac{1}{1-u^i} \right)^2, \]
which is at most
\begin{eqnarray*}
& & (n+1) \sum_{m=0}^n q^{n-m} Coef. \ u^m \ in \ \left( \prod_i \frac{1}{1-u^i} \right)^2 \\
& = & q^n (n+1) \sum_{m=0}^n \frac{1}{q^m} Coef. \ u^m \ in \ \left( \prod_i \frac{1}{1-u^i} \right)^2 \\
& \leq &  q^n (n+1) \sum_{m=0}^{\infty} \frac{1}{q^m} Coef. \ u^m \ in \ \left( \prod_i \frac{1}{1-u^i} \right)^2 \\
& = & q^n (n+1) \prod_{i \geq 1} \left( \frac{1}{1-1/q^i} \right)^2 \\
& \leq & 12(n+1)q^n.
\end{eqnarray*}

For part 3, we argue as in part 2 but at the final step use the bound (recall that $q \geq 2$) 
\[ \prod_{i \geq 1}  \left( \frac{1}{1-1/q^i} \right)^2 \leq (1+22/q).\]

For part 4, it follows from Theorem \ref{finegen} that the number of orbits divided by $q^n$ is the coefficient of $u^n$ in
\[ \left( \frac{1}{1-u} \right)^2 \prod_{i \geq 1} \left( \frac{1}{1-u^{i+1}/q^i} \right)^2.\] Since $\prod_{i \geq 1} \left( \frac{1}{1-u^{i+1}/q^i} \right)^2$
is analytic in a disk of radius larger than 1, part 4 follows from standard asymptotics (Theorem 11.7 of \cite{Od}).
\end{proof}

\section{Orthogonal groups}

Here we consider the orthogonal groups acting on the Lie algebra of the special orthogonal group.   We have several subsections corresponding to nilpotent and general orbits
in odd characteristic and characteristic two.    
The case of characteristic two is much more challenging.

\subsection{Odd characteristic: nilpotent orbits}

Over an algebraically closed field, it is well known that the nilpotent orbits are uniquely determined by their Jordan canonical form and the only restriction on
the Jordan form is that blocks of even size occur in pairs. In this section, we work over a finite field of odd characteristic. We treat only the odd dimensional orthogonal groups (the even case is similar but more involved since we need to consider the type of the form).

The number of nilpotent orbits of $O(2n+1,q)$ on its Lie algebra is equal to the number of unipotent conjugacy classes of $O(2n+1,q)$.  Lemma \ref{counter} gives a generating function for this number. This is implicit in \cite{Wall}.

\begin{lem} \label{counter} Suppose that the characteristic is odd. For $n \geq 0$, let $a_n^{\pm}$ be the number of nilpotent orbits of $O^{\pm}(n,q)$ on its Lie algebra. 
\begin{enumerate}
\item \[ 1 + \sum_{n \geq 1} (a_n^+ + a_n^-) u^n = \prod_{i \ even} \frac{1}{1-u^{2i}} \prod_{i \ odd} \frac{1+u^i}{1-u^i}.\]
\item \[ 1 + \sum_{n \geq 1} (a_n^+ - a_n^-) u^n = \prod_{i \ even} \frac{1}{1-u^{2i}}.\]
\item In particular, the number of nilpotent orbits  of $O(2n+1,q)$ on its Lie algebra is equal to

\[ \left[ \textrm{Coeff. \ } u^{2n+1} \textrm{\ in \ } \frac{1}{2} \prod_{i \ even} \frac{1}{1-u^{2i}} \prod_{i \ odd} \frac{1+u^i}{1-u^i} .\right] \]
\end{enumerate}
\end{lem}

\begin{proof} We work with unipotent conjugacy classes of $O^{\pm}(n,q)$. These correspond to partitions of $n$ with all even sizes occurring with even multiplicity, and with a choice of sign for the distinct odd sizes such that the overall sign is $\pm$. Thus

\begin{eqnarray*}
1 + \sum_{n \geq 1} (a_n^+ + a_n^-) u^n & = & \prod_{i \ even} (1+u^{2i}+u^{4i}+ \cdots) \prod_{i \ odd} (1+2u^i+2u^{2i}+2u^{3i}+\cdots) \\
& = & \prod_{i \ even} \frac{1}{1-u^{2i}} \prod_{i \ odd} \left( 1 + \frac{2u^i}{1-u^i} \right) \\
& = & \prod_{i \ even} \frac{1}{1-u^{2i}} \prod_{i \ odd} \frac{1+u^i}{1-u^i},
\end{eqnarray*} proving part 1.

Similarly, part 2 follows since
\[ 1 + \sum_{n \geq 1} (a_n^+ - a_n^-) u^n = \prod_{i \ even} (1+u^{2i}+u^{4i}+ \cdots) \prod_{i \ odd} 1 = \prod_{i \ even} \frac{1}{1-u^{2i}},\]
and part 3 is obvious from part 1 since $a_{2n+1}^+ = a_{2n+1}^-$. \end{proof}

 \subsection{Odd characteristic: general orbits}

Once we have the formula for the nilpotent case, obtaining the generating function for the total number of orbits is not too difficult. Since we are working with odd dimension, for each characteristic polynomial which could occur, write it as $x^e h(x)$ with $h(0) \ne 0$ and so $e$ is odd. Letting $a_e$ be as in Lemma \ref{counter}, the number of orbits is given by multiplying $a_e$ by what happens in the odd characteristic symplectic case with no power of $x$ occurring. This leads to the following result.

\begin{thm} \label{countemall} Suppose that the characteristic is odd. Then the number of orbits of $O(2n+1,q)$ on its Lie algebra is equal to the coefficient of $u^{2n+1}$
in \[ \frac{1}{2} \prod_{i \geq 1} \frac{1}{1-qu^{2i}} \prod_{i \ odd} (1+u^i)^2.\] 
\end{thm}

\begin{proof} By the above remarks and arguing as in the proof of Theorem \ref{sympoddgen}, the total number of orbits is the coefficient of $u^{2n+1}$ in 
\begin{equation}  B(u) \prod_{d \geq 1} (\sum_{\lambda} u^{2d |\lambda|})^{N^*(q;2d)} \prod_{d \geq 1}  (\sum_{\lambda} u^{2d |\lambda|})^{M^*(q;d)}, \end{equation}
where by Lemma \ref{counter} \[ B(u) =  \frac{1}{2} \prod_{i \ even} \frac{1}{1-u^{2i}} \prod_{i \ odd} \frac{1+u^i}{1-u^i}.\]

Again arguing as in Theorem \ref{sympoddgen} and then rewriting, this is equal to
\begin{eqnarray*}
B(u) \prod_{i \geq 1} \frac{1-u^{2i}}{1-qu^{2i}} & = & \frac{1}{2}  \prod_{i \geq 1} \frac{1-u^{2i}}{1-qu^{2i}} \prod_{i \ even} \frac{1}{1-u^{2i}}
\prod_{i \ odd} \frac{(1+u^i)^2}{1-u^{2i}} \\
& = & \frac{1}{2} \prod_{i \geq 1} \frac{1}{1-qu^{2i}} \prod_{i \ odd} (1+u^i)^2.
\end{eqnarray*}
\end{proof}

{\it Remark:} Setting $q=1$ in Theorem \ref{countemall} gives the generating function for the number of nilpotent orbits.

Theorem \ref{excel} gives excellent upper and lower bounds on the number of orbits of $O(2n+1,q)$ on its Lie algebra.

\begin{thm} \label{excel} Suppose that the characteristic is odd.
\begin{enumerate}
\item The number of orbits of $O(2n+1,q)$ on its Lie algebra is at most $3.5 q^n$.
\item The number of orbits of $O(2n+1,q)$ on its Lie algebra is at most $q^n + 7q^{n-1}$.
\item  The number of orbits of $O(2n+1,q)$ on its Lie algebra is at least $q^n$.
\end{enumerate}
\end{thm}

\begin{proof} By Theorem \ref{countemall}, the sought quantity is the coefficient of $u^{2n}$ in 
\[ \frac{1}{2u} \prod_{i \geq 1} \frac{1-u^{2i}}{1-qu^{2i}} \frac{\prod_{i \ odd} (1+u^i)^2}{\prod_i (1-u^{2i})}.\]
This is equal to
\begin{eqnarray*}
& & \textrm{Coeff. \ } u^{2n} \textrm{\ in \ } \frac{1}{4u} \prod_{i \geq 1} \frac{1-u^{2i}}{1-qu^{2i}}
\left[  \frac{\prod_{i \ odd} (1+u^i)^2}{\prod_i (1-u^{2i})} - \frac{\prod_{i \ odd} (1-u^i)^2}{\prod_i (1-u^{2i})} \right] \\
& = & \frac{1}{4} \sum_{m=0}^n  \textrm{Coeff. \ } u^{2n-2m} \textrm{\ in \ } \prod_{i \geq 1} \frac{1-u^{2i}}{1-qu^{2i}} \\
& & \cdot  \textrm{Coeff. \ } u^{2m} \textrm{\ in \ } \frac{1}{u} \left[  \frac{\prod_{i \ odd} (1+u^i)^2}{\prod_i (1-u^{2i})} -
\frac{\prod_{i \ odd} (1-u^i)^2}{\prod_i (1-u^{2i})} \right].
\end{eqnarray*}

As in the proof of Theorem \ref{glbound}, this is at most 
\[ \frac{1}{4} q^n \sum_{m=0}^n \frac{1}{q^m} \textrm{Coeff. \ } u^{2m} \textrm{\ in \ } \frac{1}{u} \left[ \frac{\prod_{i \ odd} (1+u^i)^2}{\prod_i (1-u^{2i})} -
\frac{\prod_{i \ odd} (1-u^i)^2}{\prod_i (1-u^{2i})} \right]. \] Since every term in the power series expansion of  \[ \frac{1}{u} \left[  \frac{\prod_{i \ odd} (1+u^i)^2}{\prod_i (1-u^{2i})} -
\frac{\prod_{i \ odd} (1-u^i)^2}{\prod_i (1-u^{2i})} \right] \] is non-negative (the bracketed part is of the form $f(u)-f(-u)$ where all terms of the power series expansion of $f(u)$ are non-negative), we obtain an upper bound of
\[ \frac{1}{4}q^n q^{1/2} \left[ \frac{\prod_{i \ odd} (1+1/q^{i/2})^2}{\prod_i (1-1/q^i)} -   \frac{\prod_{i \ odd} (1-1/q^{i/2})^2}{\prod_i (1-1/q^i)}  \right].\] 

Using Euler's pentagonal number theorem, which gives arbitrarily good upper and lower bounds on infinite products such as $\prod_j (1-1/q^j)$, and noting that our products over odd $i$ can be expressed in terms of such products, it is easy to see that if $q \geq 3$, then \[ \frac{1}{4} q^{1/2} \left[ \frac{\prod_{i \ odd} (1+1/q^{i/2})^2}{\prod_i (1-1/q^i)} -   \frac{\prod_{i \ odd} (1-1/q^{i/2})^2}{\prod_i (1-1/q^i)}  \right] \]  is at most $3.5$ and also at most $1+7/q$. This proves the first two parts of the theorem.

For the lower bound, note that
\[  \textrm{Coeff. \ } u^{2n+1} \textrm{\ in \ } \frac{1}{2} \prod_{i \geq 1} \frac{1}{1-qu^{2i}} \prod_{i \ odd} (1+u^i)^2 \] is at least
\[ \textrm{Coeff. \ } u^{2n} \textrm{\ in \ } \frac{1}{2} \prod_{i \geq 1} \frac{1}{1-qu^{2i}} \cdot \textrm{Coeff. \ } u \textrm{\ in \ } \prod_{i \ odd} (1+u^i)^2,\]
which is at least $q^n$.
\end{proof}

\subsection{Even characteristic: nilpotent orbits}

In this subsection, we will prove the results stated in the introduction for nilpotent orbits in characteristic $2$. 

We first deal with the algebraically closed case and the modify the results for the finite case.  Since the centralizers are not in general connected,
the classes over the algebraic closure will split in some number of classes in the finite case (always a power of $2$) and in the $-$ case
some orbits will not be present.   By \cite{LS}, the number of nilpotent orbits is independent of $q$.  

Canonical forms are given for the nilpotent classes over the algebraic closure in \cite[Chapter 5]{LS}.  

Any such element can be written as an orthogonal direct sum of indecomposable elements of the form
$W(m)$ and $W_{\ell}(n)$ subject to various constraints.  This converts the problem to a purely combinatorial problem.

\begin{enumerate}
\item For each positive integer $m$, we have an object $W(m)$ of weight $m$. We
can use any number of copies of these objects.

\item For each pair of positive integers $(\ell,n)$ with $(n+1)/2<\ell\le n$,
we have an object $W_{\ell}(n)$ of weight $n$. We can use multiple objects
of this type, but the $n$-values corresponding to them must be distinct.
Furthermore, if we use objects with indices $(\ell_i,n_i)$, where we arrange
them so that the sequence $(n_i)$ is in (strictly) decreasing order, then the
sequences $(\ell_i)$ and $(n_i-\ell_i)$ are also strictly decreasing.
\end{enumerate}

Let $d_n$ be the number of ways to choose a collection of objects subject to
these conditions with total weight $n$ and let
$$D(q)=\sum_{n=0}^{\infty} d_n q^n$$
be the corresponding generating function. As before let
$P(q)=\prod_{n=1}^{\infty} \frac{1}{1-q^n}$ be the generating function for the
partition function.   Note that $d_n$ is the number of nilpotent orbits of the orthogonal group
on $\mathfrak{so}(2n,k)$ with $k$ algebraically closed of characteristic $2$. 
 
 \begin{lem} \label{countnil2} 
$D(q)=(1-q) P(q)^2.$
\end{lem}

This shows that $(d_n)$ is the sequence A000990 from the OEIS and also the difference
sequence $d_n = a_n-a_{n-1}$ of the previous case.

\begin{proof} 
If we define the sequence $\ell'_i=n_i-\ell_i$, then we have
$0\le \ell'_i<(n_i-1)/2$ for all $i$ and the sequences $(\ell_i)$ and
$(\ell'_i)$ are strictly decreasing. Thus putting the $\ell_i$ on the top row
and the $\ell'_i$ on the bottom row, we can regard these as giving a way to fill
a $2\times r$ grid with non-negative integers (for some $r$), so that they are
strictly decreasing as we go right and as we go down they decrease by at least
2. Let $G'_r$ be the set of all such $2\times r$ grids. Conversely, if we
take such a grid $g\in G'_r$ and use its rows to define $\ell'_i$ and $\ell_i$, then $n_i=\ell_i+\ell'_i$ will be strictly decreasing, we will have
$2\ell_i>\ell_i+\ell'_i+1=n_i+1$, and the total weight $\sum_i n_i$ will be
just the sum $|g|$ of all the entries of $g$. It is routine to write down the
generating function corresponding to these pieces. For a $2\times r$ table,
there is a minimal example with top row $(r+1,r,\dots,2)$ and bottom row
$(r-1,r-2,\dots,0)$ which has entries totalling $r(r+1)$. Subtracting this
from $g$, will give a $2\times r$ grid that is nonincreasing as we go right
or down, and such grids are up to cosmetic changes plane partitions based on a
$2\times r$ grid. This gives
$$\sum_{g\in G'_r} q^{|g|}=\frac{q^{r(r+1)}}{(1-q)(1-q^2)^2\cdots (1-q^r)^2(1-q^{r+1})}.$$
Hence the full generating function for objects of the second class is
$$\sum_{r=0}^{\infty} \frac{q^{r(r+1)}}{(1-q)(1-q^2)^2\cdots (1-q^r)^2(1-q^{r+1})}.$$
Since the first class contributes an independent factor of $P(q)$ to the 
generating function, the claimed formula just reduces to the following lemma.
\end{proof} 

\begin{lem} 
$\sum_{r=0}^{\infty} \frac{q^{r(r+1)}}{(1-q)^2(1-q^2)^2\cdots (1-q^r)^2(1-q^{r+1})} = P(q)$. 
\end{lem} 

\begin{proof} View a partition $\lambda=(\lambda_1,\lambda_2,\dots)$
with $\lambda_1\ge \lambda_2 \ge \cdots$ as a Ferrers diagram with $\lambda_i$
cells in row $i$. Thus cells correspond to pairs $(i,j)$ with
$1\le j\le \lambda_i$ for all $i$. The trick is just to break up diagrams
according to the largest $r$ such that they contain the cell $(r,r+1)$. If a
diagram contains the cell $(r,r+1)$, then it must contain the entire
rectangle $[1,r]\times [1,r+1]$, which is a total of $r(r+1)$ cells. Since it
does not contain $(r+1,r+2)$, what is left after we cut out this rectangle is
a smaller diagram with at most $r$ rows and another with at most $r+1$
columns. Conversely, given any diagrams with at most $r$ rows and at most
$r+1$ columns, respectively, we can attach them to the $r\times (r+1)$
rectangle to get exactly such a diagram. The generating functions for these
pieces are standard. Having at most $r$ rows gives
$$\frac{1}{(1-q)(1-q^2)\cdots (1-q^r)}$$
and having at most $r+1$ columns is after transposing exactly the same (with
$r$ incremented of course). Thus these diagrams have generating function
$$\frac{q^{r(r+1)}}{(1-q)^2(1-q^2)^2\cdots (1-q^r)^2(1-q^{r+1})},$$
and we are done.
\end{proof} 

We now turn to the count over a finite field of characteristic $2$.   Each class defined
above splits into a certain number of classes over the finite field.   Aside from the case where
the class contains elements which are an orthogonal sum of $W(m_i)$ with all $m_i$ even,
the number of classes is the same for both types.  In this case, the classes only appear in $+$ type
(see \cite{LS}).  

In particular, we see that if $n$ is odd, then there are the same number of nilpotent classes for both types of
orthogonal Lie algebras.   

Now let us turn to the weighted version.  Each object above gets a weight of 
$2^{r+s}$ where:\\

\noindent (1) $r$ is the number of odd indices that occur among the $m$ values
corresponding to the first class of objects and are not in any of the 
intervals $[2(n_i-\ell_i),2\ell_i-1]$. \\

\noindent (2) $s$ is the number of indices $i$ such that $\ell_{i-1}+\ell_i\le n_{i-1}$,
where provided there is at least one such block we count this inequality as
holding vacuously for $i=1$.\\

By \cite[5.12]{LS},  the weighted count is the sum of the number of classes of nilpotent 
orbits for the two types of orthogonal Lie algebras.  

\begin{lem}  The generating function which enumerates the above objects is
the square of the generating function for the number of partitions.
\end{lem} 

\begin{proof} 
This count, modulo working out some details, is very similar to the symplectic case.
If we focus on just the first class of blocks and we ignore the possible
exclusions, then the generating function is just
$$\frac{\prod_{n~\rm odd} (1+q^n)}{\prod_{n=1}^{\infty} (1-q^n)}.$$
The reason for this is that for odd $n$ the total contribution is
$$\frac{1+q^n}{1-q^n}=1+2q^n+2q^{2n}+2q^{3n}+\cdots,$$
that is we get an extra factor of 2 for every odd part that occurs in our
collection of blocks (and for even $n$ there is no such factor so we just get
$\frac{1}{1-q^n}$.)

For the second class, the basic objects are again grids $G'_r$. The condition
$\ell_{i-1}+\ell_i\le n_{i-1}$ translates to $\ell'_{i-1}\ge \ell_i$, so we
define $s(g)$ to be the number of indices for which $\ell'_{i-1}\ge \ell_i$,
where by convention $\ell'_0=\infty$ so $i=1$ always counts as such an index
(if $r\ge 1$ so that $1$ is an index we allow). The excluded interval, in terms
of $\ell_i$ and $\ell'_i$ is the odd numbers in the range
$[2\ell'_i,2\ell_i-1]$, which we will write as numbers of the form $2x-1$
where $x\in [\ell'_i+1,\ell_i]$. Accordingly, for a grid $g$ define the 
excluded set to be 
$$X_g = \bigcup_{i=1}^r [\ell'_i+1,\ell_i].$$
Then the contribution of the second class of blocks is
$$E(q)=\sum_{r=0}^{\infty} \sum_{g\in G'_r} \frac{2^{s(g)} q^{|g|}}{\prod_{i\in X_g} (1+q^{2i-1})},$$
where the denominator is chosen to cancel off the extra factors of $1+q^n$
for the odd parts which are to be excluded from the count. This sum is 
intimidating but the same bag of tricks applies. Let $E_{k,n}$ be the total
contribution to $E(q)$ from all grids (of whatever size) whose first column
is $\left(\begin{matrix} n\cr k\end{matrix}\right)$. Then for $k=0$ only the
$2\times 1$ grid with just this column contributes and hence we get
$$E_{0,n}(q) = \frac{2q^n}{\prod_{i=1}^n (1+q^{2i-1})},$$
and pulling off the first column immediately gives the recurrence
$$E_{k,n}(q) = \frac{2q^{n+k}}{\prod_{i=k+1}^n (1+q^{2i-1})} + \sum_{j=0}^{k-1} \sum_{m=j+2}^{n-1} \frac{2^{\delta_{k\ge m}} q^{n+k} E_{j,m}}{\prod_{i=1+\max(k,m)}^n (1+q^{2i-1})}.$$
The first term corresponds to the single column case. Otherwise, when we peel
off the first column we get a smaller grid $g'$ and the sum breaks the
contribution up according to the first column of $g'$. The factor of $2$ in the
numerator is exactly recording where the first two columns give an extra
contribution to $s$ and the denominator is the observation that for any such
$g$ and $g'$, these are exactly the exclusions that hold for $g$ but not for
the smaller grid $g'$.

The first step in the calculation is to prove that
$$E_{k,n}(q) = \frac{2 q^{n+k} \prod_{i=n}^{n+k-1} (1-q^i)}{\prod_{i=1}^k (1-q^i) \prod_{i=1}^n (1+q^{2i-1})}$$
satisfies this recursion, hence induction shows that this formula holds for
all $k,n$. This amounts to showing that for all $k$, $n$ (with $k\le n-2$), we
have
\begin{equation*}\begin{split}
& \frac{2q^{n+k}\prod_{i=n}^{n+k-1} (1-q^i)}{\prod_{i=1}^k (1-q^i) \prod_{i=1}^n (1+q^{2i-1})} = \frac{2q^{n+k}}{\prod_{i=k+1}^n (1+q^{2i-1})}\\
& \qquad + \sum_{j=0}^{k-1} \sum_{m=j+2}^{n-1} \frac{2^{1+\delta_{k\ge m}} q^{n+k+m+j} \prod_{i=m}^{m+j-1} (1-q^i)}{\prod_{i=1}^j (1-q^i) \prod_{i=1}^m (1+q^{2i-1}) \prod_{i=1+\max(k,m)}^n (1+q^{2i-1})},
\end{split}\end{equation*}
which after cancelling common factors and clearing denominators reads
\begin{equation*}\begin{split}
& \prod_{i=n}^{n+k-1} (1-q^i) = \prod_{i=1}^k (1-q^i) \prod_{i=1}^k (1+q^{2i-1})\\
&\qquad +  \sum_{j=0}^{k-1} \sum_{m=j+2}^{n-1} 2^{\delta_{k\ge m}} q^{m+j} \prod_{i=m}^{m+j-1} (1-q^i) \prod_{i=j+1}^{k} (1-q^i)  \prod_{i=m+1}^k (1+q^{2i-1}).
\end{split}\end{equation*}
Here again we take the convention that the final product is empty (hence equal to $1$) if $m+1> k$. Interchanging the order of the sums, which is most clearly done by splitting the $m$ sum into two pieces, makes this
\begin{equation*}\begin{split}
& \prod_{i=n}^{n+k-1} (1-q^i) = \prod_{i=1}^k (1-q^i) \prod_{i=1}^k (1+q^{2i-1})\\
&\qquad +  \sum_{m=2}^k \sum_{j=0}^{m-2} 2 q^{m+j} \prod_{i=m}^{m+j-1} (1-q^i) \prod_{i=j+1}^k (1-q^i)  \prod_{i=m+1}^k (1+q^{2i-1})\\
&\qquad +  \sum_{m=k+1}^{n-1} \sum_{j=0}^{k-1} q^{m+j} \prod_{i=m}^{m+j-1} (1-q^i) \prod_{i=j+1}^k (1-q^i).
\end{split}\end{equation*}
Using the identity
$$q^j \prod_{i=m}^{m+j-1}(1-q^i) \prod_{i=j+1}^k (1-q^i) = \prod_{i=m+1}^{m+j}(1-q^i) \prod_{i=j+1}^k (1-q^i) -  \prod_{i=m+1}^{m+j-1}(1-q^i) \prod_{i=j}^k (1-q^i),$$
which holds since after cancelling off common factors it reduces to
$$q^j(1-q^m)=(1-q^{m+j}) - (1-q^j),$$
we see that both $j$ sums telescope, hence we reduce to
\begin{equation*}\begin{split}
& \prod_{i=n}^{n+k-1} (1-q^i) = \prod_{i=1}^k (1-q^i) \prod_{i=1}^k (1+q^{2i-1})\\
&\qquad +  \sum_{m=2}^k 2 q^m \prod_{i=m+1}^{2m-2} (1-q^i) \prod_{i=m-1}^k (1-q^i)  \prod_{i=m+1}^k (1+q^{2i-1})\\
&\qquad +  \sum_{m=k+1}^{n-1} q^m (1-q^k) \prod_{i=m+1}^{m+k-1} (1-q^i).
\end{split}\end{equation*}
The first sum telescopes since
\begin{equation*}\begin{split}
& 2 q^m \prod_{i=m+1}^{2m-2} (1-q^i) \prod_{i=m-1}^k (1-q^i)  \prod_{i=m+1}^k (1+q^{2i-1})\\
& = \prod_{i=m+1}^{2m} (1-q^i) \prod_{i=m+1}^k (1-q^i)(1+q^{2i-1}) -
\prod_{i=m}^{2m-2} (1-q^i) \prod_{i=m}^k (1-q^i)(1+q^{2i-1}),
\end{split}\end{equation*}
(which after cancelling common factors, and splitting $1-q^{2m}=(1-q^m)(1+q^m)$ to find an extra common factor, reads $2q^m(1-q^{m-1})=(1-q^{2m-1})(1+q^m)-(1-q^m)(1+q^{2m-1})$) and the second telescopes since
$$q^m(1-q^k)\prod_{i=m+1}^{m+k-1} (1-q^i) = \prod_{i=m+1}^{m+k} (1-q^i) - \prod_{i=m}^{m+k-1} (1-q^i).$$
Plugging these in and seeing that everything cancels completes the proof.

Having proved the formula for $E_{k,n}(q)$, we have a double sum for the 
generating function corresponding to blocks of the second type
$$E=1+\sum_{n=2}^{\infty} \sum_{k=0}^{n-2} \frac{2q^{n+k}\prod_{i=n}^{n+k-1} (1-q^i)}{\prod_{i=1}^k (1-q^i) \prod_{i=1}^n (1+q^{2i-1})}.$$
Writing
$$q^k = \frac{(1-q^{n+k})-(1-q^k)}{1-q^n},$$
we see that
\begin{equation*}\begin{split}
E & = 1 + \sum_{n=2}^{\infty} \frac{2q^n}{\prod_{i=1}^n (1+q^{2i-1})} \sum_{k=0}^{n-2} \left[\frac{\prod_{i=n+1}^{n+k} (1-q^i)}{\prod_{i=1}^k (1-q^i)}-\frac{\prod_{i=n+1}^{n+k-1} (1-q^i)}{\prod_{i=1}^{k-1} (1-q^i)}\right]\\
& = 1 + \sum_{n=2}^{\infty} \frac{2q^n \prod_{i=n+1}^{2n-2} (1-q^i)}{\prod_{i=1}^{n-2} (1-q^i) \prod_{i=1}^n (1+q^{2i-1})}\\
& = 1 + \sum_{n=2}^{\infty} \left[\frac{\prod_{i=n+1}^{2n} (1-q^i)}{\prod_{i=1}^n (1-q^i)(1+q^{2i-1})}-\frac{\prod_{i=n}^{2n-2} (1-q^i)}{\prod_{i=1}^{n-1} (1-q^i)(1+q^{2i-1})}\right]\\
& =\lim_{N\to\infty} \frac{\prod_{i=N+1}^{2N} (1-q^i)}{\prod_{i=1}^N (1-q^i)(1+q^{2i-1})}=\prod_{i=1}^{\infty} \frac{1}{(1-q^i)(1+q^{2i-1})}.
\end{split}\end{equation*}
This can be rewritten in many ways. One nice one is to first write
$$E=\prod_{i=1}^{\infty} \frac{1+q^i}{(1-q^{2i})(1+q^{2i-1})},$$
and then note that the odd index numerators exactly cancel all the factors of $1+q^{2i-1}$. Hence we are left with only the even index numerators, which we can compress to give
$$E=\prod_{i=1}^{\infty} \frac{1+q^{2i}}{1-q^{2i}}.$$
This form makes it clear that $E$ has only even powers of $q$.

Combining this with the generating function we already found for the first 
class of blocks, we get the full generating function for these two block types
with the power of 2 weights of
$$ A(q)=\frac{\prod_{n~\rm odd} (1+q^n)}{\prod_{n=1}^{\infty} (1-q^n)} \cdot \prod_{i=1}^{\infty} \frac{1+q^{2i}}{1-q^{2i}}.$$
The terms in the numerator combine to be $\prod_n (1+q^n)$, which is the
generating function for partitions into distinct pieces and since these are
well known to be equinumerous with the partitions into odd pieces, we get
$$A(q) = \prod_{n=1}^{\infty} \frac{1}{1-q^n}\cdot \prod_{i=1}^{\infty} \frac{1}{1-q^{2i}}\cdot \prod_{i=1}^{\infty} \frac{1}{1-q^{2i-1}} = \prod_{n=1}^{\infty} \frac{1}{(1-q^n)^2},$$
completing the proof.
\end{proof} 

Given our discussion before the lemma, we immediately have:

\begin{cor} The number of nilpotent $O^+(2n,q)$
orbits on $\mathfrak{so}^+(2n,q)$ plus the number of $O^-(2n,q)$ nilpotent orbits on $\mathfrak{so}^-(2n,q)$ is the number of $Sp(2n,q)$ nilpotent orbits
on $\mathfrak{sp}(2n,q)$. 
\end{cor}

The difference in the class numbers of the two orthogonal Lie algebras is very easy to compute.  As noted above, it is only the classes corresponding to the orthogonal direct
sum of $W(m_i)$ with all $m_i$  even.  The centralizers of elements in these classes are all connected and so the  difference in the number of such classes in 
 is $0$ for $n$ odd and is  $p(n/2)$ for $n$ even.  These classes do not appear in $-$ type and so:

\begin{lem}  The number of nilpotent $O^+(2n,q)$ orbits on $\mathfrak{so}^+(2n,q)$ minus the number of nilpotent $O^-(2n,q)$ orbits on $\mathfrak{so}^-(2n,q)$ is $0$
if $n$ is odd and $p(n/2)$ if $n$ is even.  
\end{lem}

Let $d_n$ denote the number of nilpotent orbits in the symplectic case
and $e_n^{\epsilon}$, the number of $O^{\epsilon}(2n,q)$ orbits on nilpotents in the corresponding orthogonal Lie algebra. 
The result above gives 
$e_n^+=e_n^- = d_n/2$ if $n$ is odd and $e_n^+ - e_n^- = p(n/2)$ if $n$ is even.  Thus, we have:

\begin{thm} \label{soeven} 
\begin{enumerate} 
\item $d_n = e_n^+ + e_n^-.$; 
\item if $n$ is odd, $e_n^+=e_n^- = d_n/2$.
\item If $n$ is even  $e_n^{\epsilon} = (1/2)(d_n + \epsilon \cdot p(n/2))$.
\end{enumerate}
\end{thm}

If we consider $SO$ orbits, the only quantity that changes is the number of orbits when  $n$ is even and we are working with 
$\mathfrak{so}^+(2n,q)$.  The number of orbits increases by $p(n/2)$.

\subsection{Even characteristic: general orbits}

We next consider semisimple orbits.   The characteristic polynomial of any element is $f^2$ for a monic polynomial $f$ of degree $n$.  Writing 
$f^2 = \prod_{i=1}^r f_i^{2e_i}$ where the $f_i$ are distinct monic irreducible polynomials,  we see that there are $2^{r-1}$ orbits in each orthogonal Lie
algebra.  This can be seen by writing the space $V$ as an orthogonal direct sum of $V_i, 1 \le i \le r$ where the characteristic polynomial is $f_i^{2e_i}$.  
Then $V_i$ is a nondegenerate space of type $\epsilon$.  This corresponds to the signed polynomial $f$ where $f_i$ is assigned the sign $\epsilon_i$.   
A semisimple element with a given characteristic polynomial and choice of signs  is contained in precisely one of the types of Lie algebras given by the product of the signs.  

We first prove a combinatorial lemma that is valid over any finite field (but does not have the same interpretation in odd characteristic).  By a signed polynomial,
we mean a polynomial together with a choice of sign for each distinct irreducible factor.  

\begin{thm} For $n \geq 1$, the number of degree $n$ monic signed polynomials over the field $F_q$ is $(n+1)q^n - (n-1)q^{n-1}$.  
\end{thm}

\begin{proof} We need to count each monic polynomial with coefficients in $F_q$ by the weight $2^r$ where $r$ is the number of distinct irreducible factors of the polynomial. Clearly this is equal to

\[ \textrm{Coeff. \ } u^n \textrm{\ in }\prod_{\phi} (1+2u^{deg(\phi)}+2u^{2 deg(\phi)} \cdots), \] where the product is over all monic irreducible polynomials $\phi$
with coefficients in $F_q$. This is equal to 

\begin{eqnarray*}
\textrm{Coeff. \ }u^n \textrm{\ in} \prod_{\phi} \left(1+\frac{2u^{deg(\phi)}}{1-u^{deg(\phi)}}\right) & = & \textrm{Coeff. \ }u^n \textrm{\ in} \prod_{\phi} \left(\frac{1+u^{deg(\phi)}}{1-u^{deg(\phi)}} \right) \\
& = & \textrm{Coeff. \ } u^n \textrm{ \ in} \prod_{\phi} \frac{1-u^{2deg(\phi)}}{(1-u^{deg(\phi)})^2} \\
& = &  \textrm{ Coeff. \ } u^n \textrm{ \ in} \frac{1-u^2q}{(1-uq)^2},
\end{eqnarray*}  where the final step used Lemma \ref{necklace} twice.

Clearly this is equal to

\begin{eqnarray*}
 & & \textrm{Coeff.\ } u^n \textrm{\ in} \frac{1}{(1-uq)^2} - \textrm{Coeff. \ } u^n \textrm{\ in} \frac{u^2q}{(1-uq)^2} \\
& = & q^n \cdot \textrm{Coeff.\ }u^n \textrm{\ in} \frac{1}{(1-u)^2} - q \cdot \textrm{Coeff. \ } u^{n-2} \textrm{\ in} \frac{1}{(1-uq)^2} \\
& = & q^n (n+1) - q q^{n-2} \cdot \textrm{Coeff. \ } u^{n-2} \textrm{\ in} \frac{1}{(1-u)^2} \\
& = & (n+1)q^n - (n-1)q^{n-1}.
\end{eqnarray*}
\end{proof}

Note that this implies that if $r$ is the number of distinct irreducible factors of a degree $n$ polynomial, then for large $q$, the expected value of $2^r/(n+1)$ is close to $1$. 
 
As noted, 

\begin{thm} For $n \geq 1$, the number of semisimple orbits  of $O^{\epsilon}(2n,q)$ in $\mathfrak{so}^{\epsilon}(2n,q)$ is $ ((n+1)q^n - (n-1)q^{n-1})/2$.
\end{thm}

We now consider general orbits.  The simplest quantity to compute is the sum of the two class numbers.

\begin{thm} \label{sum}  Let $q$ be even. The  number of orbits of $O^+(2n,q)$ on $\mathfrak{so}^+(2n,q)$ plus the number of orbits of $O^-(2n,q)$ on $\mathfrak{so}^-(2n,q)$   is
the number of orbits of $Sp(2n,q)$ on $\mathfrak{sp}(2n,q)$.
\end{thm}

\begin{proof}  Let $f(x)$ be a monic polynomial of degree $n$ with $f= \prod_{i=1}^r f_i^{e_i}$ with the $f_i$ distinct irreducible polynomials.
Let $n_i = e_i \deg(f_i)$.  
Let $h^{\epsilon}(i)$ denote the number of  $O^{\epsilon}(2e_i,q)$ orbits  of nilpotent elements on $\mathfrak{so}^{\epsilon}(2n_i,q)$  (note
that this number does not change if we replace $q$ by a any power of $q$).  

We will prove a stronger result where the equality holds for each $f$.  If we fix a choice of signs for each irreducible factor, we see that
the number of nilpotent classes in the centralizer of a semisimple element with signed characteristic polynomial $f^2$ is 
$\prod_{i=1}^r h^{\epsilon_i}(e_i)$.    Then we add these over all possible choices of signs we obtain $\prod_{i=1}^r (h^+(e_i) + h^-(e_i)) =
\prod_{i=1}^r d(e_i)$ where $d(e_i)$ is the number of nilpotent orbits in $\mathfrak{sp}(2e_i,q)$.
 The latter quantity is precisely the number of nilpotent classes in the centralizer of a semisimple element 
in $\mathfrak{sp}(2n,q)$ with characteristic polynomial $f^2$.   Summing over all possible $f$ gives the result. 
\end{proof} 

Similar reasoning gives a nice formula for the difference in the number of classes.

\begin{thm}  \label{diff}   Let $\Delta(n)$ be the difference between the number of (orthogonal group) orbits on $\mathfrak{so}^+(2n,q)$ and the number of orbits on $\mathfrak{so}^-(2n,q)$.
Then $\Delta(n)=0$ if $n$ is odd and $\Delta(n)$ is the number of orbits of $GL(n/2,q)$ on $\mathfrak{gl}(n/2,q)$ if $n$ is even.    
\end{thm}

\begin{proof}  Let $f(x)$ be a monic polynomial of degree $n$ with $f= \prod_{i=1}^r f_i^{e_i}$ with the $f_i$ distinct irreducible polynomials.  Define $n_i$ and
$h^{\epsilon}(i)$ as in the previous proof.

We will prove a stronger result where the equality holds for each $f$.  If we fix a choice of signs for each irreducible factor, we see that
the number of nilpotent classes in the centralizer of a semisimple element with characteristic polynomial $f^2$ and fixed choice of signs for each $i$ is 
$\prod_{i=1}^r h^{\epsilon}(e_i)$.    Then we add these over all possible choices of signs to obtain $\prod_{i=1}^r (h^+(e_i) - h^-(e_i))$.
If at least one of the $e_i$ is odd, then we see the difference is $0$.  In particular if $n$ is odd, the difference is $0$.

Thus, we only need consider $f$ with $f=g^2$ (so the characteristic polynomial is a $4$th power).  In particular all $e_i$ are even and so
we see that the difference is $\prod_{i=1}^r (h^+(e_i) - h^-(e_i)) = \prod_{i=1}^r p(e_i/2)$.  Note that this is the same as computing the number of 
nilpotent classes commuting with a semisimple element in $\mathfrak{gl}(n/2,q)$ with characteristic polynomial $g$.  Summing over all possible $g$
gives the number of $GL(n/2,q)$ orbits on $\mathfrak{gl}(n/2,q)$.   This completes the proof. 
\end{proof} 

We thus have generating functions for both the sum and difference and so for each term as well and we see that the number of classes for $\mathfrak{so}^{\epsilon}(2n,q)$  is 
very close to $(1/2)$ the number of orbits on $\mathfrak{sp}(2n,q)$ (exactly if $n$ is odd).   We leave it to the reader to use the bounds for $\mathfrak{sp}$ to write down
similar bounds.  We note a fairly good lower bound comes from just considering semisimple elements.

To end this section, we note that in characteristic $2$, $Sp(2n,q)$ acts on each type of orthogonal Lie algebra (but not on the orthogonal groups). By \cite[Theorem 71]{Kap}, 
the $Sp$ orbits are determined by their rational canonical form and each possible rational canonical form is just the double of a rational
canonical form in $\mathfrak{gl}(n,q)$.     Note that Kaplansky assumes the field is algebraically closed.  However, we note the elementary fact:

\begin{lem}  Let $k$ be an algebraically closed field of characteristic $2$.  Let $X=Sp(2n,k)$ acting on $\mathfrak{sp}(2n,k)$.  If $x \in \mathfrak{sp}(2n,k)$, then
$C_X(x)$ is connected.
\end{lem}

\begin{proof}  Note that the centralizer of a semisimple element in the Lie algebra is a direct product of symplectic groups.  So by induction, the problem reduces to the
nilpotent case.  This result is noted in \cite{LS}.  
\end{proof} 

In particular by the Lang-Steinberg theorem, the intersection of each $X$-orbit with $\mathfrak{sp}(2n,q)$ is an $X(q)$-orbit.  

This immediately gives the following result.

\begin{thm} \label{Sp on so}   Let $q$ be a power of $2$.   Let $G=Sp(2n,q)$ act on $\mathfrak{so}^{\epsilon}(2n,q)$.
\begin{enumerate}
\item The number of nilpotent orbits of $G$ on $\mathfrak{so}^{\epsilon}(2n,q)$ is $p(n)$.  
\item  The number of orbits of $G$ on $\mathfrak{so}^{\epsilon}(2n,q)$ is the number of orbits of $GL(n,q)$ on $\mathfrak{gl}(n,q)$.
\end{enumerate}
\end{thm} 

 \section{Commuting Varieties}
 
 In this section, we discuss the structure of the commuting subvariety of $G \times \mathfrak{g}$. where $G$ is an algebraic group
 and $\mathfrak{g} $ is the Lie algebra of $G$ (i.e. the variety of pairs $(g,x)$ with $x^g=x$).  
 We show that our earlier results determine the dimension of this subvariety and the number of components of maximal dimension.
 
 So let $k$ be an algebraically closed field of characteristic $p > 0$.  We first note that if $G=GL(n,k)$, then since $G$ is an open dense subvariety
 of $\mathfrak{gl}(n,k)$, it follows by Motzkin-Taussky \cite{MT} that the variety is irreducible of dimension $n^2 + n$ (and is the closure of 
 $(T \times \mathfrak{t})^G$ where $\mathfrak{t}$ is the Lie algebra of a maximal torus).   
 
 It is known by Levy \cite{Le} that if $\mathfrak{g}$ is a reductive Lie algebra in very good characteristic, then the variety of commuting pairs in
 $\mathfrak{g} \times \mathfrak{g}$ is an irreducible variety of dimension equal to $\dim \mathfrak{g} + \mathrm{rank}(G)$ (and the description is
 as in the $GL$ case).     This was proved by Richardson \cite{Ri} in characteristic zero (and the proof in positive characteristic is very similar --- it
 was also shown in \cite{RG} that the characteristic zero results implies the very good characteristic case).  
 
 This is consistent with our results on the number of orbits.   This does not immediately imply the result for $G \times \mathfrak{g}$.  However, the same proof 
 with minor modifications does hold.
 
Finally, we consider the most interesting case for $p=2$ and $G=Sp(2n,k)$, $O(2n,k)$ or $SO(2n,k)$.  We note that Roman  \cite{Ro} has studied the commuting variety 
of $\mathfrak{g} \times \mathfrak{g}$ where the dimension is bigger than in odd characteristic and is not irreducible in the case of the orthogonal group.

By the Lang-Weil bound, it follows that for $G=Sp(2n,k)$ or $G=O(2n,k)$, the dimension of the commuting variety is $\dim G + n$ and the number of components 
of maximal dimension  is $n+1$.   We will now describe the $n+1$ components of maximal dimension.  

First consider $G=Sp(2n,k)$.  Let $V$ be the natural module for $G$ and write $V = V_1 \perp \ldots \perp V_n$ where $\dim V_i=2$. 
Let $e_i \ne 0$ be a nilpotent root element acting on each $V_s$  and acting trivially if $s \ne i$. 
Fix $j$ with $ 0 \le j \le n$.   Let $S_j \subset \mathfrak{g}$ consist  of all elements preserving each $V_i$ and acting as $a_iI + d_ie_i$ where
$d_i = 1$ if $i \le j$ and $d_i=0$ otherwise.  Moreover, we assume that the $a_i$ are distinct.   Then the centralizer in $Sp(2n,k)$ also preserves
each $V_i$ and either induces $Sp(2,k)$ on $V_i$ if $i > j$ or a $1$-dimensional unipotent subgroup  if $i \le  j$.  Let $X_i$ be the closure of 
$(C_i \times S_i)^G$ where $C_i$ is the centralizer of $S_i$ (note all elements of $S_i$ have the same centralizer).  It is straightforward (by computing generic fibers)
that $\dim X_i = \dim G + n$ and $X_i$ is irreducible.  Moreover clearly $X_i \ne  X_s$ for $i \ne s$ (and there are open disjoint subsets in each).  
Thus, the $X_i, 0 \le i \le n$ are the distinct components of maximal dimension.  

Next consider $G=O(2n,k)$.   Let $V$ and $V_i$ be as above.  Let $T =T_1 \times \ldots \times T_n$ where $T_i= SO(V_i)$  is a maximal torus.
Let $\mathfrak{t}$ be the Lie algebra of $T$.  
Fix an involution $w_j \in G$ with $w_j$ inverting  $T_1, \ldots T_j$ and centralizing $T_i$ for $i > j$.   Let $S$ be the subset of the Lie algebra $\mathfrak{t}$ 
so that the $V_i$ are the distinct eigenspaces.   For each $j$ with $0 \le j \le n$, consider the closure $X_i$ of $(w_jT \times S)^G$.  Again, we see that
$X_i$ is irreducible of dimension $\dim G + n$ and the $X_i$ are distinct and so they are precisely the components of maximal dimension.

Note that if $n$ is odd, then precisely $(n+1)/2$ of the components are contained in $SO(2n,k) \times \mathfrak{g}$ and $(n+1)/2$ are not.
If $n$ is even, then precisely $(n+2)/2$ of the components are contained in $SO(2n,k) \times \mathfrak{g}$.  We summarize:

\begin{thm}  Let $k$ be an algebraically closed field of characteristic $2$.  Let $G$ be one of the algebraic group $Sp(2n,k)$, $O(2n,k)$ or $SO(2n,k)$.
Let $\mathfrak{g}$ be the Lie algebra of $G$ and let $X_G$ the commuting subvariety of $G \times \mathfrak{g}$. 
\begin{enumerate}
\item $\dim X_G = \dim G  + n$.
\item  If $G=Sp(2n,k)$ or $O(2n,k)$,  $X_G$ has exactly $n+1$ components of maximal dimension.  
\item  If $G = SO(2n,k)$, then $X_G$ has $(n+1)/2$ components of maximal dimension if $n$ is odd and $(n+2)/2$ components
of maximal dimension if $n$ is even.
\end{enumerate} 
\end{thm}

One can also consider the commuting subvariety $X$ of $Sp(2n,k) \times \mathfrak{so}(2n,k)$.   Given the remarks in the previous section, we
see that $\dim X = \dim Sp(2n,k) + n$ and there is a unique component of maximal dimension.  In fact, in this case, we can show that 
$X$ is irreducible.  

\begin{thm} Let $k$ be an algebraically closed field of characteristic $2$.  Let $G=Sp(2n,k)$.  Then the commuting subvariety $X$ of $Sp(2n,k) \times \mathfrak{so}(2n,k)$
is irreducible of  $\dim G + n$.  Let $T$ be maximal torus of $SO(2n,k)$ and let $\mathfrak{t}$ be its Lie algebra.  Let $Y = (T \times \mathfrak{t})^G$.  Then 
$X$ is the closure of $Y$.  
\end{thm} 

\begin{proof} Note that $\dim Y = \dim X$ and so the closure of $Y$ is a component of $X$ of maximal dimension.   We show that any commuting pair
$(g, w) \in X$ is contained in the closure of $Y$.  We induct on $n$.  If $n=1$, the result is clear (in this case $X= G \times \mathfrak{so}(2n,k)$).  

If $g$ and $w$ both preserve a proper nondegenerate space, the result follows by induction.  So assume this is not the case.   In particular, this implies 
$w = \lambda I + w'$ where $w'$ is nilpotent.   So it reduces to the case that $w$ is nilpotent.  Suppose first that $w$ has more than two Jordan blocks
(and so at least four).   Then there exists $h \in C_G(w)$ with $h$ having at least four distinct eigenvalues.  Then $(h,w) \in \bar{Y}$ for any such $h$
since $h$ fixes a canonical proper nondegenerate space which must also be $w$-invariant.  Since the set of elements in $C_G(w)$ having at least four 
distinct eigenvalues is an open subset of $C_G(w)$ and $C_G(w)$ is irreducible, this set is also dense in $C_G(w)$.  Thus an open dense subset
of $C_G(w) \times \{ w \} \subset \bar{Y}$ and so  $C_G(w) \times \{ w \} \subset \bar{Y}$.   

So it remains to consider the case that $w$ has exactly two Jordan blocks (of size $n$).  We note that $\dim C_G(w) = 3n$.  Note that if $s \in \mathfrak{t}$, 
has $n$ distinct eigenvalues, then $\dim C_G(s)=3n$.   Note that the projection map $p$ from $\bar{Y}$ to $\mathfrak{so}(2n,k)$ is surjective (restrict
to the Levi subalgebra $\mathfrak{gl}(n,k)$) and  thus the generic fiber of $p$ has dimension $3n$. Since $C_G(w)$ 
is connected of dimension $3n$, it follows that the full fiber of $p:X \rightarrow \mathfrak{so}(2n,k)$ is contained in $\bar{Y}$, whence
$C_G(w) \times \{ w \}$ is contained in $\bar{Y}$ and so $X = \bar{Y}$.
\end{proof}


\begin{thebibliography}{AAA}

\bibitem{DH} Diaconis, P. and Howes, M., Random sampling of contingency tables and partitions: two practical examples of the Burnside process, {\it Statistics and Computing}, Vol. 35 (2026), Article 181.

\bibitem{DZ} Diaconis, P. and Zhong, C., Counting the number of group orbtis by marrying the Burnside process with importance sampling, {\it Adv. Appl. Math.} {\bf 172} (2026), article 102955, 24 pages.

\bibitem{FG} Fulman, J. and Guralnick, R., Bounds on the number and sizes of conjugacy classes in finite Chevalley groups with applications
to derangements, {\it Trans. Amer. Math. Soc.} {\bf 364} (2012), 3023-3070.

\bibitem{He} Hesselink, W. H. Nilpotency in classical groups over a field of characteristic 2,  {\it Math. Z.}  {\bf 166} (1979), 165--181. 

\bibitem{Kap} Kaplansky, I., Linear Algebra and Geometry,  A Second Course,  Reprint of the 1974 revised edition.
Dover Publications, Inc., Mineola, NY, 2003. 

\bibitem{Le} Levy, P.,  Commuting varieties of Lie algebras over fields of prime characteristic, {\it J. Algebra} {\bf  250} (2002),  473--484.  
 
\bibitem{LS} Liebeck, M. and Seitz, G., Unipotent and nilpotent classes in simple algebraic groups and Lie algebras, American Mathematical Society, 2012

\bibitem{MR} Maslen, D. and Rockmore, D., Separation of variables and the computation of Fourier transforms on finite groups, I., {\it J.
Amer. Math. Soc.} {\bf 10} (1997), 169-214.

\bibitem{MT} Motzkin, T. and Taussky, O.,  Pairs of matrices with property L, II, {\it Trans. Amer. Math. Soc.}  {\bf 80}  (1955), 387--401. 

\bibitem{Od} Odlyzko, A., Asymptotic enumeration methods, in {\it Handbook of Combinatorics Volume 2}, 1995, 1063-1230.

\bibitem{Ri}  Richardson, R. W.,  Commuting varieties of semisimple Lie algebras and algebraic groups  {\it Compositio Math.}  {\bf 38}  (1979), 311--327. 

\bibitem{Ro} Roman, V.,  Commuting varieties of Lie algebras in bad characteristic,  preprint.  

\bibitem{RG} Roman, V. and Guralnick, R.,   The variety of nilpotent piairs $(A,B)$ with $[A,B] = \lambda I$,  {\it Lin. Algebra Appl} {\bf  736} (2026),  214--224.



\bibitem{St} Stanley, R., Enumerative combinatorics, Volume 1, Second edition, Cambridge University Press, 2012.

\bibitem{TZ}  Taussky, O. and Zassenhaus, H.,  On the similarity transformation between a matrix and its transpose, 
{\it Pacific J. Math.} {\bf  9}  (1959), 893--896.

\bibitem{Wall} Wall, G. E., On the conjugacy classes in the unitary, symplectic and orthogonal groups, 
 {\it J. Austral. Math. Soc.}  3 (1963), 1--62.
 
 \bibitem{Xu}    Xue, T.,  Nilpotent orbits in classical Lie algebras over fields of characteristic $2$ and the Springer correspondence,  
 {\it Represent. Theory}  {\bf 13}  (2009), 371---390.

\end{thebibliography}
\end{document}